\documentclass[11pt]{article}
\usepackage{epic,latexsym,amssymb}
\usepackage{color}
\usepackage{tikz}
\usepackage{comment}
\usepackage{lineno}
\usepackage{pstricks,pst-node,pst-plot}
\usetikzlibrary{patterns}
\usepackage{geometry,graphicx}

\usetikzlibrary{calc}

\usepackage{amsmath}

\newtheorem{theorem}{Theorem}
\newtheorem{lemma}{Lemma}

\newtheorem{problem}{Problem}
\newtheorem{corollary}{Corollary}
\newtheorem{observation}{Observation}
\newtheorem{proposition}{Proposition}

\newcommand{\QED}{$\Box$}

\newcommand{\cH}{\mathcal{H}}

\newcommand{\cT}{\mathcal{T}}

\newcommand{\modo}{{\rm mod \,}}

\newcommand{\planar}{{\rm planar}}
\newcommand{\outplanar}{{\rm outerplanar}}

\newcommand{\barH}{\overline{H}}

\newcommand{\alphaF}{\alpha_{\mathcal{F}}}

\newcommand{\cG}{{\cal G}}
\newcommand{\cub}{{\rm cubic}}
\newcommand{\bipartite}{{\rm bipartite}}

\newcommand{\proof}{\noindent\textbf{Proof. }}
\newcommand{\2}{ \vspace{0.2cm} }
\newcommand{\1}{ \vspace{0.1cm} }

\let\oldenumerate\enumerate
\renewcommand{\enumerate}{
  \oldenumerate
  \setlength{\itemsep}{0pt}
  \setlength{\parskip}{0pt}
  \setlength{\parsep}{0pt}
}

\newtheorem{staticthm}{Theorem}
\newcounter{staticthmctr}
\begin{document}


\title{Independence, induced subgraphs, and \\
domination in $K_{1,r}$-free graphs}

\author{$^a$Yair Caro, $^{b,c}$Randy Davila\footnote{Corresponding author}, $^d$Michael A. Henning and $^{b,e}$Ryan Pepper\\
\\
$^a$ Department of Mathematics \\
University of Haifa--Oranim\\
Tivon 36006, Israel \\
\small \tt Email: yacaro@kvgeva.org.il \\
\\[-0.35em]
$^b$ First Principles \\
100 University Ave \\
Toronto, ON M5J 1V6, Canada \\
\small \tt Email: randy@firstprinciples.org \\
\\[-0.35em]
$^c$ Department of Computational Applied \\ Mathematics \& Operations Research\\
Rice University \\
Houston, TX 77005, USA \\
\small \tt Email: rrd6@rice.edu \\
\\[-0.35em]
$^d$ Department of Mathematics and Applied Mathematics \\
University of Johannesburg \\
Auckland Park, 2006 South Africa\\
\small \tt Email: mahenning@uj.ac.za \\
\\[-0.35em]
$^e$
Department of Mathematics and Statistics \\
University of Houston--Downtown \\
Houston, TX 77002, USA \\
\small \tt Email: pepperr@uhd.edu
}

\date{}
\maketitle

\begin{abstract}
Let $G$ be a graph and $\mathcal{F}$ a family of graphs. Define $\alpha_{\mathcal{F}}(G)$ as the maximum order of any induced subgraph of $G$ that belongs to the family $\mathcal{F}$. For the family $\mathcal{F}$ of graphs with \emph{chromatic number} at most~$k$, we prove that if $G$ is $K_{1,r}$-free, then $\alpha_{\mathcal{F}}(G) \le (r-1)k\gamma(G)$, where $\gamma(G)$ is the \emph{domination number}. When $\mathcal{F}$ is the family of empty graphs, this bound simplifies to $\alpha(G) \le 2\gamma(G)$ for $K_{1,3}$-free (claw-free) graphs, where $\alpha(G)$ is the \emph{independence number} of $G$. For $d$-regular graphs, this is further refined to the bound $\alpha(G) \le 2\left(\frac{d+1}{d+2}\right)\gamma(G)$, which is tight for $d \in \{2, 3, 4\}$. For the \emph{$k$-independence number} $\alpha_k(G)$, we prove that if $G$ is $K_{1,r}$-free with order $n$ and minimum degree $\delta \ge k+1$,
\[
\alpha_k(G) \le \left( \frac{(r-1)(k+1)}{\delta - k + (r-1)(k+1)} \right) n,
\]
and this bound is sharp for all parameters. Finally, using a local Ramsey-domination lemma, we extend this framework to edge-hereditary graph families, showing that for $K_{1,r}$-free graphs, we have $\alpha_{\mathcal{F}}(G) \le (r(\mathcal{F}^*,K_r) - 1)\gamma(G)$, where $\mathcal{F}^*$ is the set of graphs not in $\mathcal{F}$. Specializing to $K_q$-free graphs, we show $\alpha_{\mathcal{F}}(G) \le (r(K_q,K_r) - 1)\gamma(G)$.
\end{abstract}

{\small \textbf{Keywords:}  Claw-free graphs, domination number, independence number,  $K_{1,r}$-free graphs} \\
\indent {\small \textbf{AMS subject classification:} 05C69}

\section{Introduction}
\label{sec:introduction}

The study of independence and domination in graphs is a central topic in graph theory, with numerous applications and significant structural insights. Let $G$ be a finite, simple graph with vertex set $V(G)$. A set $S \subseteq V(G)$ is \emph{independent} if no two vertices in $S$ are adjacent, and the \emph{independence number} $\alpha(G)$ is the maximum cardinality among all independent sets in $G$. A set $D \subseteq V(G)$ is a \emph{dominating set} if every vertex in $V(G) \setminus D$ has at least one neighbor in $D$, where two vertices are neighbors if they are adjacent. The \emph{domination number} $\gamma(G)$ is the minimum cardinality among all dominating sets in $G$. The relationship between $\alpha(G)$ and $\gamma(G)$, particularly in restricted graph classes, has been extensively studied; see, for instance, the domination monograph~\cite{HaHeHe-23}, results on forbidden substructures~\cite{Brause-2022}, independence in claw-free and $4$-regular graphs~\cite{Kang-2014}, and conjectures on independence and domination generated by the computer program \emph{TxGraffiti}~\cite{CaDaHePe2022a,Davila-2026-TxGraffiti}. Generalized domination and independence parameters, including $k$-domination and $j$-independence, have also been studied in this broader context~\cite{Pepper-2010,HansbergPepper-2013}.

The independence number $\alpha(G)$ can be viewed as a special case of a more general problem concerning induced subgraphs. Let $\mathcal{F}$ be a family of graphs, and define $\mathcal{F}(G)$ as the set of induced subgraphs of $G$ that belong to $\mathcal{F}$. Let $\alpha_{\mathcal{F}}(G)$ denote the maximum order of any element in $\mathcal{F}(G)$. This perspective motivates the problem of bounding $\alpha_{\mathcal{F}}(G)$ in terms of $\gamma(G)$, particularly for graphs that exclude certain subgraphs.


\subsection{Graph theory notation and terminology}
\label{sec:notation}

For notation and graph theory terminology, we generally follow~\cite{HaHeHe-23}. Specifically, let $G$ be a graph with vertex set $V(G)$ and edge set $E(G)$, and of order~$n(G) = |V(G)|$ and size $m(G) = |E(G)|$. A \emph{neighbor} of a vertex $v$ in $G$ is a vertex $u$ adjacent to $v$, that is, $uv \in E(G)$. The \emph{open neighborhood} $N_G(v)$ of a vertex $v$ in $G$ is the set of neighbors of $v$, while the \emph{closed neighborhood} of $v$ is the set $N_G[v] = \{v\} \cup N(v)$. We denote the \emph{degree} of $v$ in $G$ by $\deg_G(v) = |N_G(v)|$. For a set $S \subseteq V(G)$, its \emph{open neighborhood} is the set $N_G(S) = \cup_{v \in S} N_G(v)$, and its \emph{closed neighborhood} is the set $N_G[S] = N_G(S) \cup S$. The \emph{chromatic number} of $G$, written $\chi(G)$, is the minimum number of colors to assign to the vertices of $G$ so that no two adjacent vertices share the same color.

For a set $S \subseteq V(G)$, the subgraph induced by $S$ is denoted by $G[S]$. Further, the subgraph of $G$ obtained from $G$ by deleting all vertices in $S$ and all edges incident with vertices in $S$ is denoted by $G - S$; that is, $G-S = G[V(G)\setminus S]$. If $S = \{v\}$, then we also denote $G - S$ simply by $G - v$. If $F$ is a graph, then an $F$-\emph{component} of $G$ is a component isomorphic to~$F$. We will also use the notation $xG$ to denote the disjoint union of $x \ge 1$ copies of the graph $G$. The notation $xG + H$ denotes the graph obtained by attaching each vertex of the graph $H$ to every vertex of $xG$. For an integer $k \ge 1$, we will also use the standard notation $[k] = \{1, \dots, k\}$.

For any two graphs $G$ and $H$, the \emph{generalized Ramsey number} $r(G,H)$ is the smallest integer $p$ for which a red subgraph isomorphic to $G$ or a blue subgraph isomorphic to $H$ appears in any red-blue edge coloring of the complete graph $K_p$. By ``red subgraph'' we mean a subgraph whose edges are all colored red, and similarly for the phrase ``blue subgraph.'' The special case $r(K_m,K_n)$ is the classical Ramsey number, also denoted $r(m,n)$ in the literature. When $\mathcal{F}$ is a family of graphs, then the Ramsey number $r(G,\mathcal{F})$ is the smallest integer $p$ for which a red subgraph isomorphic to $G$ or a blue subgraph isomorphic to a graph that belongs to the family $\mathcal{F}$ appears in any red-blue edge coloring of the complete graph $K_p$; the notation $r(\mathcal{F},G)$ is defined analogously, with the family in the red color and the graph $G$ in the blue color. 

We denote the \emph{path}, \emph{cycle}, and \emph{complete graph} on $n$ vertices by $P_n$, $C_n$, and $K_n$, respectively, and we denote the \emph{complete bipartite graph} with partite sets of cardinality~$n$ and $m$ by $K_{n,m}$. A \emph{triangle} in $G$ is a subgraph isomorphic to $K_3$. A graph is \emph{$F$-free} if it does not contain $F$ as an induced subgraph. In particular, if $F = K_{1,3}$, then the graph is \emph{claw-free}, while if $F = K_4-e$, then the graph is \emph{diamond-free}. An excellent survey of claw-free graphs has been written by Flandrin, Faudree, and Ryj\'{a}\v{c}ek~\cite{claw_free_survey}. Chudnovsky and Seymour attracted considerable interest in claw-free graphs due to their excellent series of papers in \textit{Journal of Combinatorial Theory} on this topic (see, for example, their paper~\cite{claw_free_JCTB}). For $r \ge 3$ we define the class $\cG_{r}$ of graphs as follows, where we remark that in the special case when $r = 3$, the class $\cG_{r}$ is the well-studied class of claw-free graphs.

\paragraph{The class $\cG_{r}$.} For $r \ge 3$, let $\cG_{r}$ be the class of all $K_{1,r}$-free graphs.

\section{Main results}

In this section we state our main results. Our first result establishes general bounds on $\alpha_{\mathcal{F}}(G)$, where $\mathcal{F}$ is the family of graphs with chromatic number at most~$k$.

\begin{theorem}
\label{thm:general_intro_bound}
For $r \ge 3$, if $G \in \cG_{r}$ and $\mathcal{F}$ is the family of graphs with chromatic number at most $k$, then
\[
\alpha_{\mathcal{F}}(G) \le (r-1)k\gamma(G),
\]
and this bound is sharp.
\end{theorem}

We note that if $\mathcal{F}$ is the family of graphs with chromatic number at most~$1$, then $\mathcal{F}$ consists of edgeless graphs, and in this case $\alpha_{\mathcal{F}}(G) = \alpha(G)$, yielding the bound $\alpha(G) \le (r-1)\gamma(G)$. Moreover for claw-free graphs (that is, when $r = 3$), this further simplifies to $\alpha(G) \le 2\gamma(G)$. The family $\mathcal{F}$ of graphs with chromatic number at most~$k$ where $k \in \{2,3,4\}$ corresponding to bipartite ($k=2$), outerplanar ($k=3$), and planar graphs ($k=4$), respectively, and in these cases we establish the sharp bound of $\alpha_{\mathcal{F}}(G) \le (r-1)k\gamma(G)$.

We next extend our analysis to \emph{$k$-independence numbers} in $K_{1, r}$-free graphs. A set $S \subseteq V(G)$ is \emph{$k$-independent} if the subgraph of $G$ induced by $S$ has maximum degree at most~$k$. The \emph{$k$-independence number}, $\alpha_k(G)$, is the maximum cardinality among all $k$-independent sets in $G$. We prove the following sharp bound for $K_{1,r}$-free graphs, which generalizes results in~\cite{Faudree-92,LiVi-90}.

\begin{theorem}
\label{thm:general_intro_kindep}
For $k \ge 0$ and $r \ge 3$, if $G \in \cG_r$ has order $n$ and minimum degree $\delta \ge k+1$, then
\[
\alpha_k(G) \le \left( \frac{(r-1)(k+1)}{\delta - k + (r-1)(k+1)} \right) n,
\]
and this bound is sharp for all parameters involved, including for every $n \equiv 0 \pmod{(r-1)(k+1) + \delta - k}$.
\end{theorem}

For induced subgraphs with unbounded chromatic number, a local Ramsey-domination lemma yields the following.

\begin{theorem}
\label{thm:general_intro_ramsey_kq}
For integers $r, q \ge 3$, if $G \in \cG_{r}$ and if $\mathcal{F}$ is the family of $K_q$-free graphs, then
\[
\alphaF(G) \le \left(r(K_q,K_r) - 1\right)\gamma(G),
\]
and this bound is sharp.
\end{theorem}

Extending this framework, we provide bounds for induced subgraphs belonging to edge-hereditary graph families. A family of graphs $\mathcal{F}$ is \emph{edge}-\emph{hereditary} if for every graph $H \in \mathcal{F}$, every subgraph of $H$ also belongs to $\mathcal{F}$. Let $\mathcal{F}^*$ denote the set of graphs not in $\mathcal{F}$. We prove the following general result.

\begin{theorem}
\label{thm:general_intro_ramsey}
For $r \ge 3$, if $G \in \cG_{r}$ and $\mathcal{F}$ is an edge-hereditary family of graphs and $\mathcal{F}^*$ denotes the set of graphs not in $\mathcal{F}$, then
\[
\alphaF(G) \le \left(r(\mathcal{F}^*,K_r) - 1\right) \gamma(G),
\]
and this bound is sharp.
\end{theorem}

The remainder of this paper is organized as follows.
Section~\ref{sec:known} presents known results foundational to our analysis. In Section~\ref{sec:k-chromatic}, we study the interplay between maximum $k$-chromatic induced subgraphs and domination in $K_{1,r}$-free graphs, with applications to independence. Section~\ref{sec:k-independence} establishes a computable bound on the $k$-independence number in $K_{1,r}$-free graphs. Section~\ref{sec:ramsey} provides bounds on the maximum cardinality of $K_q$-free induced subgraphs in $K_{1,r}$-free graphs using Ramsey numbers and domination. Finally, Section~\ref{sec:conclusion} offers concluding remarks and open problems.

\section{Known results}
\label{sec:known}

In this section, we recall results on independence and domination that will be used to analyze the bounds we prove. Every vertex in a dominating set in a graph $G$ dominates itself and at most~$\Delta(G)$ other vertices, yielding the following trivial lower bound on the domination number.

\begin{observation}
\label{obs:bound-dom}
If $G$ is a graph of order~$n$, then $\gamma(G) \ge \frac{n}{\Delta(G) + 1}$.
\end{observation}

As a consequence of Brook's Coloring Theorem, if $G \ne K_n$ is a connected graph of order~$n$ with maximum degree~$\Delta \ge 3$, then the chromatic number of $G$ is at most~$\Delta$, that is, $\chi(G) \le \Delta$. Alternatively, viewing a $\Delta$-coloring of $G$ as a partitioning of its vertices into $\Delta$ independent sets, called color classes, we infer by the Pigeonhole Principle that $G$ contains an independent set of cardinality at least~$n/\Delta$, implying that $\alpha(G) \ge n/\Delta$. We state this formally as follows.

\begin{observation}
\label{obs:lower-bound-alpha-1}
If $G \ne K_{n}$ is a connected graph of order~$n$ and $\Delta(G) \ge 3$, then $\alpha(G) \ge \frac{n}{\Delta(G)}$.
\end{observation}

Several authors have given the following upper bound on the independence number of a claw-free graph, which is the $r=3$ and $k=0$ instance of Theorem~\ref{thm:general_intro_kindep}.
\begin{theorem}{\rm (\cite{Faudree-92, LiVi-90})}
\label{thm:upper-bound-alpha}
If $G$ is a claw-free graph of order~$n$ with minimum degree~$\delta$, then
\[
\alpha(G) \le \left( \frac{2}{\delta + 2} \right) n.
\]
\end{theorem}

The following property of connected, claw-free, cubic graphs is established in~\cite{HeLo-12}.

\begin{lemma}{\rm (\cite{HeLo-12})}
\label{lem:known}
If $G \ne K_4$ is a connected, claw-free, cubic graph of order $n$, then the vertex set $V(G)$ can be uniquely partitioned into sets, each of which induces a triangle or a diamond in $G$.
\end{lemma}

By Lemma~\ref{lem:known}, the vertex set $V(G)$ of connected, claw-free, cubic graph $G \ne K_4$  can be uniquely partitioned into sets, each of which induces a triangle or a diamond in $G$. Following the notation introduced in~\cite{HeLo-12}, we refer to such a partition as a \emph{triangle}-\emph{diamond partition} of $G$, abbreviated $\Delta$-D-partition. We call every triangle and diamond induced by a set in our $\Delta$-D-partition a \emph{unit} of the partition. A unit that is a triangle is called a \emph{triangle-unit}, and a unit that is a diamond is called a \emph{diamond-unit}. (We note that a triangle-unit is a triangle that does not belong to a diamond.)

\section{On largest $k$-chromatic induced subgraphs of $K_{1,r}$-free graphs}
\label{sec:k-chromatic}

This section explores the relationship between the order of maximum induced subgraphs satisfying certain structural properties and the domination number of $K_{1,r}$-free graphs with applications to independence. Specifically, for $r \ge 3$ given a graph $G \in \cG_{r}$ and a family $\mathcal{F}$ of graphs, let $\mathcal{F}(G)$ denote the set of induced subgraphs of $G$ that belong to $\mathcal{F}$. We define $\alpha_{\mathcal{F}}(G)$ as the maximum order among all induced subgraphs of $G$ in $\mathcal{F}$. The general problem we consider is the following:

\begin{problem}
\label{problem:alpha-constants}
For each $r \ge 3$, determine or estimate the best possible constants $C_{\mathcal{F}}$ \emph{(}which depend only on~$\mathcal{F}$\emph{)}, such that $\alpha_{\mathcal{F}}(G) \le C_{\mathcal{F}} \gamma(G)$ for all graphs $G \in \cG_{r}$. These constants are given by
\[
C_{\mathcal{F}} = \sup_{G \in \cG_{r}}\; \frac{\alpha_{\mathcal{F}}(G)}{\gamma(G)}.
\]
\end{problem}

This problem generalizes earlier results on bounding the independence number of graphs in terms of the domination number. In this setting, the independence number $\alpha(G)$ corresponds to $\alpha_{\mathcal{F}}(G)$, when $\mathcal{F}$ is the family of empty graphs, notably a family of graphs with chromatic number one. Indeed, as the following theorem shows, we may solve the Problem~\ref{problem:alpha-constants} for graphs with bounded chromatic number. We next prove Theorem~\ref{thm:general_intro_bound}. Recall its statement.

\smallskip
\noindent \textbf{Theorem~\ref{thm:general_intro_bound}} \emph{For $r \ge 3$, if $G \in \cG_{r}$ and $\mathcal{F}$ is the family of graphs with chromatic number at most $k$, then
\[
\alpha_{\mathcal{F}}(G) \le (r-1)k\gamma(G),
\]
and this bound is sharp.
}

\noindent
\proof For $r \ge 3$, let $G \in \cG_{r}$ and let $\mathcal{F}$ be the family of graphs with chromatic number at most $k$. Let $D$ be a $\gamma$-set of $G$, and let $H$ be a maximum induced subgraph of $G$ such that $H \in \mathcal{F}(G)$. Let
\[
S = D \cap V(H) \quad \mbox{and} \quad B = D \setminus S.
\]

We first bound how many vertices of $H$ can be dominated by a single vertex of $D$. If a vertex $x \in B$ is adjacent to at least $(r-1)k+1$ vertices of $H$, then the subgraph of $H$ induced by these neighbors has chromatic number at most~$k$, and therefore contains an independent set of cardinality at least~$r$. Together with $x$, this independent set induces a forbidden $K_{1,r}$ in $G$. Hence every vertex of $B$ is adjacent to at most $(r-1)k$ vertices of $H$.

Now let $x \in S$. If $x$ is adjacent in $H$ to at least $(r-1)(k-1)+1$ vertices, and $Q$ is the subgraph of $H$ induced by these neighbors, then either $\alpha(Q) \ge r$, which again gives a forbidden $K_{1,r}$ centered at $x$, or $\alpha(Q) \le r-1$. In the latter case,
\[
\chi(Q) \ge \frac{|V(Q)|}{\alpha(Q)} > k-1,
\]
and so $\chi(Q) \ge k$. Since $x$ is adjacent to every vertex of $Q$, the induced subgraph $H[V(Q) \cup \{x\}]$ has chromatic number at least~$k+1$, a contradiction. Thus every vertex of $S$ is adjacent in $H$ to at most $(r-1)(k-1)$ vertices.

Since $D$ dominates $G$, every vertex of $H$ lies in the closed neighborhood of some vertex of $D$. Therefore
\[
\begin{array}{lcl}
|V(H)| & \le & (r-1)k|B| + \big((r-1)(k-1)+1\big)|S| \1 \\
& = & (r-1)k|D| - (r-2)|S| \1 \\
& \le & (r-1)k\gamma(G).
\end{array}
\]
This proves the desired upper bound. Moreover, equality in this bound can occur only when $S=\emptyset$, that is, when the chosen $\gamma$-set $D$ is disjoint from $V(H)$.

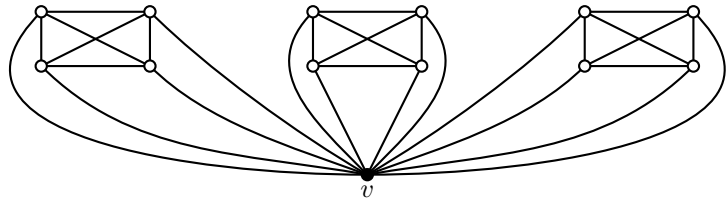
\begin{figure}[htb]
\begin{center}
\begin{tikzpicture}[scale=.85,style=thick,x=.85cm,y=.85cm]
\def\vr{2.5pt} 
\path (0,2) coordinate (u1);
\path (0,3) coordinate (u2);
\path (2,2) coordinate (u3);
\path (2,3) coordinate (u4);
\path (5,2) coordinate (v1);
\path (5,3) coordinate (v2);
\path (7,2) coordinate (v3);
\path (7,3) coordinate (v4);
\path (10,2) coordinate (w1);
\path (10,3) coordinate (w2);
\path (12,2) coordinate (w3);
\path (12,3) coordinate (w4);
\path (6,0) coordinate (v);
\draw (u1)--(u2)--(u3)--(u4)--(u1)--(u3);
\draw (u2)--(u4);
\draw (v1)--(v2)--(v3)--(v4)--(v1)--(v3);
\draw (v2)--(v4);
\draw (w1)--(w2)--(w3)--(w4)--(w1)--(w3);
\draw (w2)--(w4);
\draw (v) to[out=180,in=225, distance=2.75cm] (u2);
\draw (v) to[out=170,in=-45, distance=2cm] (u1);
\draw (v) to[out=160,in=-45, distance=1.25cm] (u3);
\draw (v) to[out=150,in=-45, distance=1.25cm] (u4);
\draw (v1)--(v)--(v3);
\draw (v) to[out=135,in=225, distance=1.25cm] (v2);
\draw (v) to[out=45,in=-45, distance=1.25cm] (v4);
\draw (v) to[out=30,in=225, distance=1.25cm] (w2);
\draw (v) to[out=20,in=225, distance=1.25cm] (w1);
\draw (v) to[out=10,in=225, distance=2cm] (w3);
\draw (v) to[out=0,in=-45, distance=2.75cm] (w4);
\draw (u1) [fill=white] circle (\vr);
\draw (u2) [fill=white] circle (\vr);
\draw (u3) [fill=white] circle (\vr);
\draw (u4) [fill=white] circle (\vr);
\draw (v1) [fill=white] circle (\vr);
\draw (v2) [fill=white] circle (\vr);
\draw (v3) [fill=white] circle (\vr);
\draw (v4) [fill=white] circle (\vr);
\draw (w1) [fill=white] circle (\vr);
\draw (w2) [fill=white] circle (\vr);
\draw (w3) [fill=white] circle (\vr);
\draw (w4) [fill=white] circle (\vr);
\draw (v) [fill=black] circle (\vr);
\draw[anchor = north] (v) node {{\small $v$}};
\end{tikzpicture}
\end{center}
\vskip -0.35cm
\caption{The $K_{1,4}$-free graph $G = 3K_4 + K_1$ in the proof of Theorem~\ref{thm:general_intro_bound}}
\label{fig:yair-easy-construction}
\end{figure}

To see that the upper bound is sharp, consider the following construction. Let $G = (r-1)K_k + K_1$, and let $v$ be the vertex of $G$ that is adjacent to all vertices in $(r-1)K_k$. In the special case when $r = 4$ and $k = 4$, the graph $G = 3K_4 + K_1$ is illustrated in Figure~\ref{fig:yair-easy-construction}, where the dominating vertex~$v$ is indicated by the shaded vertex. Every independent set in $G$ either consists of the vertex~$v$ or contains at most one vertex from each of the complete graphs $K_k$ in the $r-1$ components of $G - v$, implying that $G$ is $K_{1,r}$-free. Since $v$ is a dominating vertex of $G$, we note that $\gamma(G) = 1$. The induced subgraph $H = G - v = (r-1)K_k$ of $G$ is $k$-colorable, and so $H \in \mathcal{F}(G)$, implying that $\alpha_{\mathcal{F}}(G) \ge |V(H)| = (r-1)k = (r-1)k\gamma(G)$. By the upper bound just proved, $\alpha_{\mathcal{F}}(G) = (r-1)k\gamma(G)$. The upper bound in the statement of the theorem is therefore tight.~\QED

Theorem~\ref{thm:general_intro_bound} provides a unified framework for bounding several invariants in $K_{1,r}$-free graphs. Specifically, we can bound the maximum cardinality of induced subgraphs with restricted chromatic numbers. For example, let $\bipartite(G)$, $\outplanar(G)$, and $\planar(G)$ denote the maximum cardinality of an induced bipartite, outerplanar, and planar subgraph of $G$, respectively. These invariants, along with sharp constructions, are summarized in Table~\ref{tab:k1r_bounds}.

\begin{table}[htb]
\centering
\[
\begin{array}{|c|l|l|}
\hline
\text{Chromatic number } & \mathcal{F} & \text{General bounds for $K_{1,r}$-free graphs} \\
\hline
\chi \le 1 & \text{Empty Graphs} & \alpha(G) = \alpha_{\mathcal{F}}(G) \le (r-1)\gamma(G) \\ \hline
\chi \le 2 & \text{Bipartite Graphs} & \bipartite(G) = \alpha_{\mathcal{F}}(G) \le 2(r-1)\gamma(G) \\ \hline
\chi \le 3 & \text{Outerplanar Graphs} & \outplanar(G)=\alpha_{\mathcal{F}}(G)\le 3(r-1)\gamma(G) \\ \hline
\chi \le 4 & \text{Planar Graphs} & \planar(G)=\alpha_{\mathcal{F}}(G)\le 4(r-1)\gamma(G) \\
\hline
\end{array}
\]
\caption{Bounds on invariants in $K_{1,r}$-free graphs}
\label{tab:k1r_bounds}
\end{table}

The bounds in Table~\ref{tab:k1r_bounds} are sharp, as demonstrated by the following constructions for $K_{1,r}$-free graphs:

\begin{itemize}
    \item \textbf{Planar graphs:} Take $r-1$ disjoint copies of $K_4$, which are planar, and add a vertex $v$ adjacent to all vertices in these copies. The resulting graph is $K_{1,r}$-free, has $\gamma(G) = 1$, and removing $v$ leaves $(r-1)K_4$, and so $\planar(G) = 4(r-1)\gamma(G)$.

    \item \textbf{Outerplanar graphs:} Take $r-1$ disjoint copies of $K_3$, which are outerplanar, and add a vertex $v$ adjacent to all vertices in these copies. The resulting graph is $K_{1,r}$-free, has $\gamma(G) = 1$, and removing $v$ leaves $(r-1)K_3$, and so $\outplanar(G) = 3(r-1)\gamma(G)$.

    \item \textbf{Bipartite graphs:} Take the cycle $C_{2(r-1)}$, which is bipartite, and add a vertex $v$ adjacent to all vertices of the cycle, forming the wheel graph $W_{2(r-1)}$. The resulting graph is $K_{1,r}$-free, has $\gamma(G) = 1$, and removing $v$ leaves $C_{2(r-1)}$, and so $\bipartite(G) = 2(r-1)\gamma(G)$.
\end{itemize}

\subsection{Improvements for independence bounds in $d$-regular and claw-free graphs}

In this subsection, we refine upper bounds on the independence number $\alpha(G)$ in terms of the domination number $\gamma(G)$ for claw-free graphs, with a particular focus on $d$-regular graphs. This investigation begins with the following consequence of Theorem~\ref{thm:general_intro_bound}.

\begin{corollary}
\label{cor:main}
If $G$ is a claw-free graph, then $\alpha(G) \le 2 \gamma(G)$, and this bound is sharp.
\end{corollary}

While Corollary~\ref{cor:main} provides a general bound for claw-free graphs, improvements can be obtained in the case of $d$-regular claw-free graphs. As a consequence of Observation~\ref{obs:bound-dom} and Theorem~\ref{thm:upper-bound-alpha}, we arrive at the following sharper result.

\begin{theorem}{\rm (\cite{Faudree-92, LiVi-90})}
\label{thm:alpha-versus-dom}
For $d \ge 2$, if $G$ is a claw-free, $d$-regular graph, then
\[
\alpha(G) \le 2\left( \frac{d + 1}{d + 2} \right) \gamma(G).
\]
\end{theorem}

This result raises a natural question: for a given value of $d$, is the upper bound in Theorem~\ref{thm:alpha-versus-dom} achievable? In this section, we address this question and demonstrate that the bound is sharp for $d = 2$, $d = 3$, and $d = 4$. Specifically, we construct examples in Sections~\ref{sec:2regular},~\ref{sec:3regular}, and~\ref{sec:4regular} that achieve equality in Theorem~\ref{thm:alpha-versus-dom} for these values of $d$. However, for $d \ge 5$, the question of whether the upper bound remains achievable remains open and is posed as an open problem in the concluding section.

\subsubsection{$2$-regular graphs}
\label{sec:2regular}

When $d = 2$, Theorem~\ref{thm:alpha-versus-dom} simplifies to $\alpha(G) \le \frac{3}{2} \gamma(G)$. If $G$ is a cycle $C_n$, then $\alpha(G) = \lfloor \frac{1}{2}n \rfloor$ and $\gamma(G) = \lceil \frac{1}{3}n \rceil$. We therefore readily infer that equality holds in Theorem~\ref{thm:alpha-versus-dom} in the case when $d = 2$ when $G$ is a cycle of length congruent to zero modulo~$6$.

\begin{proposition}
\label{prop:2regular}
If $G$ is a $2$-regular graph of order~$n$, then $\alpha(G) \le \frac{3}{2} \gamma(G)$, with equality if and only if $G$ is a cycle $C_n$ where $n \equiv 0 \, (\modo \, 6)$.
\end{proposition}

\subsubsection{$3$-regular graphs}
\label{sec:3regular}

By Observation~\ref{obs:lower-bound-alpha-1}, if $G \ne K_4$ is a connected cubic graph of order $n$, then $\alpha(G) \ge \frac{1}{3}n$. Suppose further that $G \ne K_4$ is a connected, claw-free, diamond-free cubic graph of order $n$. In this case, every unit in the $\Delta$-D-partition of $G$ is a triangle-unit. Thus, $n = 3t$ where $t$ is the number of triangle-units in $G$. Moreover, the $\Delta$-D-partition of $G$ partitions $V(G)$ into $t$ sets, each of which induce a copy of $K_3$. Let $I$ be a maximum independent set of $G$, and so $|I| = \alpha(G)$. The set $I$ contains at most one vertex from every triangle-unit in our $\Delta$-D-partition, implying that $\alpha(G) = |I| \le t = \frac{1}{3}n$. As observed earlier, $\alpha(G) \ge \frac{1}{3}n$. Consequently, $\alpha(G) = \frac{1}{3}n$. This yields the following observation.

\begin{observation}
\label{obs:alpha-2}
If $G \ne K_4$ is a connected, claw-free, diamond-free cubic graph of order $n$, then $\alpha(G) = \frac{1}{3}n$.
\end{observation}

We characterize next the connected, claw-free, cubic graphs that have equal independence and domination numbers.  For this purpose, we shall need the following result on the domination number of a graph in which every vertex belongs to a triangle.

\begin{theorem}{\rm (\cite{BaHe-22a})}
\label{thm:dom-bound-cubic}
If every vertex of a graph $G$ of order~$n$ belongs to a triangle, then $\gamma(G) \le \frac{1}{3}n$.
\end{theorem}

Let $\cT_\cub$ be the family  of connected, claw-free, diamond-free cubic graphs defined as follows. For $k \ge 1$ an integer, let $F_{2k}$ be the connected cubic graph constructed as follows. Take $2k$ disjoint copies $T_1, T_2, \ldots, T_{2k}$ of a triangle, where $V(T_i) = \{x_i,y_i,z_i\}$ for $i \in [2k]$. Let
\[
\begin{array}{lcl}
E_a & = & \{x_{2i-1}x_{2i} \colon i \in [k] \} \1 \\
E_b & = & \{y_{2i-1}y_{2i} \colon i \in [k] \} \1 \\
E_c & = & \{z_{2i}z_{2i+1} \colon i \in [k] \},
\end{array}
\]
where addition is taken modulo~$2k$ (and so, $z_1 = z_{2k+1}$). Let $F_{2k}$ be obtained from the disjoint union of these $2k$ triangles by adding the edges $E_a \cup E_b \cup E_c$. The resulting graph $F_{2k}$ we call a \emph{triangle}-\emph{necklace} with $2k$ triangles. Let $\cT_\cub = \{ F_{2k} \colon k \ge 1\}$.

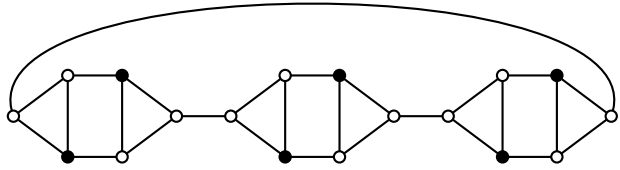
\begin{figure}[htb]
\begin{center}
\begin{tikzpicture}[scale=.85,style=thick,x=.85cm,y=.85cm]
\def\vr{2.5pt} 
\path (2,0.75) coordinate (a);
\path (3,0) coordinate (b);
\path (3,1.5) coordinate (c);
\path (4,0) coordinate (d);
\path (4,1.5) coordinate (e);
\path (5,0.75) coordinate (f);
\path (6,0.75) coordinate (g);
\path (7,0) coordinate (h);
\path (7,1.5) coordinate (i);
\path (8,0) coordinate (j);
\path (8,1.5) coordinate (k);
\path (9,0.75) coordinate (l);
\path (10,0.75) coordinate (m);
\path (11,0) coordinate (n);
\path (11,1.5) coordinate (o);
\path (12,0) coordinate (p);
\path (12,1.5) coordinate (q);
\path (13,0.75) coordinate (r);
\draw (b)--(a)--(c)--(b)--(d)--(e)--(f)--(d);
\draw (c)--(e);
\draw (f)--(g);
\draw (i)--(g)--(h)--(i)--(k)--(l)--(j)--(k);
\draw (n)--(m)--(o)--(n)--(p)--(r)--(q)--(p);
\draw (h)--(j);
\draw (l)--(m);
\draw (o)--(q);
%
\draw (a) to[out=110,in=70, distance=2.5cm] (r);
\draw (a) [fill=white] circle (\vr);
\draw (b) [fill=black] circle (\vr);
\draw (c) [fill=white] circle (\vr);
\draw (d) [fill=white] circle (\vr);
\draw (e) [fill=black] circle (\vr);
\draw (f) [fill=white] circle (\vr);
\draw (g) [fill=white] circle (\vr);
\draw (h) [fill=black] circle (\vr);
\draw (i) [fill=white] circle (\vr);
\draw (j) [fill=white] circle (\vr);
\draw (k) [fill=black] circle (\vr);
\draw (l) [fill=white] circle (\vr);
\draw (m) [fill=white] circle (\vr);
\draw (n) [fill=black] circle (\vr);
\draw (o) [fill=white] circle (\vr);
\draw (p) [fill=white] circle (\vr);
\draw (q) [fill=black] circle (\vr);
\draw (r) [fill=white] circle (\vr);
\end{tikzpicture}
\end{center}
\vskip -0.35cm
\caption{A set that is both an $\alpha$-set and a $\gamma$-set in a triangle-necklace $F_6$} \label{fig:Tneck1}
\end{figure}

Suppose that $G \in \cT_\cub$. Thus, $G$ is a triangle-necklace $F_{2k}$ for some $k \ge 1$, and so $G$ has order~$n = 6k$ and contains~$2k$ vertex disjoint triangles. We note that $G$ is diamond-free, and so Observation~\ref{obs:alpha-2} we have $\alpha(G) = \frac{1}{3}n$. Every dominating set of $G$ contains at least two vertices from every pair of triangle-units that are joined by two edges. Since $G$ contains $k$ such pairs of triangle-units, we have that $\gamma(G) \ge 2k$. Conversely, selecting one vertex from each of the $2k$ triangle-units of $G$, produces a dominating set of $G$, and so $\gamma(G) \le 2k$. Consequently, $\gamma(G) = 2k = \frac{1}{3}n = \alpha(G)$. For example, the shaded vertices in the triangle-necklace $F_6$ of order~$n = 18$ shown in Figure~\ref{fig:Tneck1} is both an $\alpha$-set and a $\gamma$-set of $F_6$ (of cardinality~$6 = \frac{1}{3}n$).

The connected, claw-free, cubic graphs achieving equality in Theorem~\ref{thm:dom-bound-cubic} are precisely the graphs in the family~$\cT_\cub$.

\begin{theorem}{\rm (\cite{BaHe-22a})}
\label{thm:bound-dom-claw-free-cubic}
If $G$ is a connected, claw-free, cubic graph of order~$n$, then $\gamma(G) \le \frac{1}{3}n$, with equality if and only if $G \in \cT_\cub$.
\end{theorem}

If $G = K_4$, then $\gamma(G) = \alpha(G) = 1$. Hence as a consequence of Observation~\ref{obs:lower-bound-alpha-1} and Theorem~\ref{thm:bound-dom-claw-free-cubic} we infer the following result.

\begin{theorem}
\label{claw-free-cubic}
If $G$ is a connected, claw-free, cubic graph of order~$n$, then $\gamma(G) \le \alpha(G)$, with equality if and only if $G \in \cT_\cub \cup \{K_4\}$.
\end{theorem}

By Observation~\ref{obs:bound-dom}, we infer that if $G$ is a cubic graph of order~$n$, then $\gamma(G) \ge \frac{1}{4}n$. By Observation~\ref{obs:alpha-2}, if $G \ne K_4$ is a connected, claw-free, diamond-free cubic graph of order $n$, then $\alpha(G) = \frac{1}{3}n$. Moreover we note that $\gamma(K_4) = \alpha(K_4) = 1$. As a consequence of these observations, we infer the following upper bound on the independence number of a claw-free, diamond-free, cubic graph in terms of its domination number.

\begin{theorem}
\label{thm:alpha-versus-dom-cubic}
If $G$ is a claw-free, diamond-free cubic graph, then $\alpha(G) \le \frac{4}{3} \gamma(G)$.
\end{theorem}

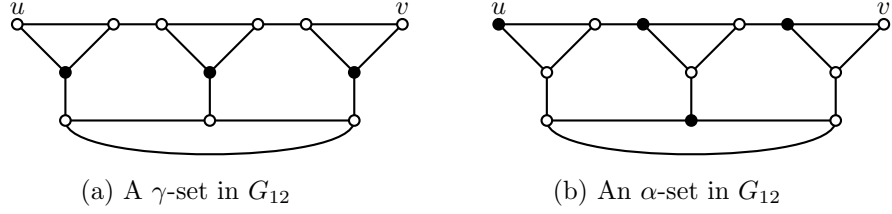
\begin{figure}[htb]
\begin{center}
\begin{tikzpicture}[scale=.8,style=thick,x=0.8cm,y=0.8cm]
\def\vr{2.5pt} 
\path (0,2) coordinate (u1);
\path (2,2) coordinate (u2);
\path (1,1) coordinate (u3);
\path (1,0) coordinate (u4);
\path (3,2) coordinate (v1);
\path (5,2) coordinate (v2);
\path (4,1) coordinate (v3);
\path (4,0) coordinate (v4);
\path (6,2) coordinate (w1);
\path (8,2) coordinate (w2);
\path (7,1) coordinate (w3);
\path (7,0) coordinate (w4);
\draw (u1)--(u2)--(u3)--(u1);
\draw (v1)--(v2)--(v3)--(v1);
\draw (w1)--(w2)--(w3)--(w1);
\draw (u3)--(u4)--(v4)--(w4);
\draw (v3)--(v4);
\draw (w3)--(w4);
\draw (u2)--(v1);
\draw (v2)--(w1);
\draw (u4) to[out=270,in=270, distance=0.75cm] (w4);
\draw (u1) [fill=white] circle (\vr);
\draw (u2) [fill=white] circle (\vr);
\draw (u3) [fill=black] circle (\vr);
\draw (u4) [fill=white] circle (\vr);
\draw (v1) [fill=white] circle (\vr);
\draw (v2) [fill=white] circle (\vr);
\draw (v3) [fill=black] circle (\vr);
\draw (v4) [fill=white] circle (\vr);
\draw (w1) [fill=white] circle (\vr);
\draw (w2) [fill=white] circle (\vr);
\draw (w3) [fill=black] circle (\vr);
\draw (w4) [fill=white] circle (\vr);
\draw[anchor = south] (u1) node {{\small $u$}};
\draw[anchor = south] (w2) node {{\small $v$}};
\draw (3.5,-1.5) node {{\small (a) A $\gamma$-set in $G_{12}$}};
\path (10,2) coordinate (x1);
\path (12,2) coordinate (x2);
\path (11,1) coordinate (x3);
\path (11,0) coordinate (x4);
\path (13,2) coordinate (y1);
\path (15,2) coordinate (y2);
\path (14,1) coordinate (y3);
\path (14,0) coordinate (y4);
\path (16,2) coordinate (z1);
\path (18,2) coordinate (z2);
\path (17,1) coordinate (z3);
\path (17,0) coordinate (z4);
\draw (x1)--(x2)--(x3)--(x1);
\draw (y1)--(y2)--(y3)--(y1);
\draw (z1)--(z2)--(z3)--(z1);
\draw (x3)--(x4)--(y4)--(z4);
\draw (y3)--(y4);
\draw (z3)--(z4);
\draw (x2)--(y1);
\draw (y2)--(z1);
\draw (x4) to[out=270,in=270, distance=0.75cm] (z4);
\draw (x1) [fill=black] circle (\vr);
\draw (x2) [fill=white] circle (\vr);
\draw (x3) [fill=white] circle (\vr);
\draw (x4) [fill=white] circle (\vr);
\draw (y1) [fill=black] circle (\vr);
\draw (y2) [fill=white] circle (\vr);
\draw (y3) [fill=white] circle (\vr);
\draw (y4) [fill=black] circle (\vr);
\draw (z1) [fill=black] circle (\vr);
\draw (z2) [fill=white] circle (\vr);
\draw (z3) [fill=white] circle (\vr);
\draw (z4) [fill=white] circle (\vr);
\draw[anchor = south] (x1) node {{\small $u$}};
\draw[anchor = south] (z2) node {{\small $v$}};
\draw (13.5,-1.5) node {{\small (b) An $\alpha$-set in $G_{12}$}};
\end{tikzpicture}
\end{center}
\vskip -0.5 cm
\caption{The claw-free, subcubic graph $G_{12}$}
\label{fig:G12}
\end{figure}

We show next that the bound in Theorem~\ref{thm:alpha-versus-dom-cubic} is achieved by an infinite family of claw-free, diamond-free cubic graphs. Let $G_{12}$ be the claw-free, subcubic graph shown in Figure~\ref{fig:G12}. The graph $G_{12}$ satisfies $\gamma(G_{12}) = 3$ and $\alpha(G_{12}) = 4$, where the shaded vertices in Figure~\ref{fig:G12}(a) and~\ref{fig:G12}(b) are examples of a $\gamma$-set and an $\alpha$-set, respectively, of $G_{12}$.

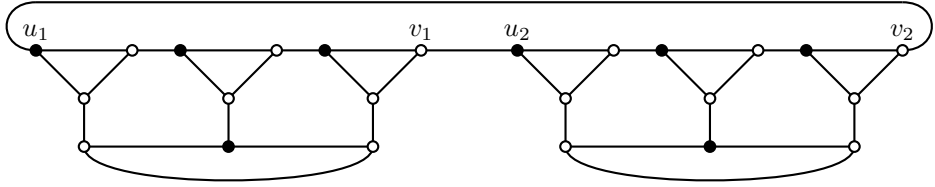
\begin{figure}[htb]
\begin{center}
\begin{tikzpicture}[scale=.8,style=thick,x=0.8cm,y=0.8cm]
\def\vr{2.5pt} 
\path (0,3) coordinate (u0);
\path (0,2) coordinate (u1);
\path (2,2) coordinate (u2);
\path (1,1) coordinate (u3);
\path (1,0) coordinate (u4);
\path (3,2) coordinate (v1);
\path (5,2) coordinate (v2);
\path (4,1) coordinate (v3);
\path (4,0) coordinate (v4);
\path (6,2) coordinate (w1);
\path (8,2) coordinate (w2);
\path (7,1) coordinate (w3);
\path (7,0) coordinate (w4);
\path (10,2) coordinate (x1);
\path (12,2) coordinate (x2);
\path (11,1) coordinate (x3);
\path (11,0) coordinate (x4);
\path (13,2) coordinate (y1);
\path (15,2) coordinate (y2);
\path (14,1) coordinate (y3);
\path (14,0) coordinate (y4);
\path (18,3) coordinate (z0);
\path (16,2) coordinate (z1);
\path (18,2) coordinate (z2);
\path (17,1) coordinate (z3);
\path (17,0) coordinate (z4);
\draw (u1)--(u2)--(u3)--(u1);
\draw (v1)--(v2)--(v3)--(v1);
\draw (w1)--(w2)--(w3)--(w1);
\draw (u3)--(u4)--(v4)--(w4);
\draw (v3)--(v4);
\draw (w3)--(w4);
\draw (u2)--(v1);
\draw (v2)--(w1);
\draw (u4) to[out=270,in=270, distance=0.75cm] (w4);
\draw (x1)--(x2)--(x3)--(x1);
\draw (y1)--(y2)--(y3)--(y1);
\draw (z1)--(z2)--(z3)--(z1);
\draw (x3)--(x4)--(y4)--(z4);
\draw (y3)--(y4);
\draw (z3)--(z4);
\draw (x2)--(y1);
\draw (y2)--(z1);
\draw (x4) to[out=270,in=270, distance=0.75cm] (z4);
\draw (w2)--(x1);
\draw (z2) to[out=0,in=0, distance=0.65cm] (z0);
\draw (u1) to[out=180,in=180, distance=0.65cm] (u0);
\draw (u0)--(z0);
\draw (u1) [fill=black] circle (\vr);
\draw (u2) [fill=white] circle (\vr);
\draw (u3) [fill=white] circle (\vr);
\draw (u4) [fill=white] circle (\vr);
\draw (v1) [fill=black] circle (\vr);
\draw (v2) [fill=white] circle (\vr);
\draw (v3) [fill=white] circle (\vr);
\draw (v4) [fill=black] circle (\vr);
\draw (w1) [fill=black] circle (\vr);
\draw (w2) [fill=white] circle (\vr);
\draw (w3) [fill=white] circle (\vr);
\draw (w4) [fill=white] circle (\vr);
\draw (x1) [fill=black] circle (\vr);
\draw (x2) [fill=white] circle (\vr);
\draw (x3) [fill=white] circle (\vr);
\draw (x4) [fill=white] circle (\vr);
\draw (y1) [fill=black] circle (\vr);
\draw (y2) [fill=white] circle (\vr);
\draw (y3) [fill=white] circle (\vr);
\draw (y4) [fill=black] circle (\vr);
\draw (z1) [fill=black] circle (\vr);
\draw (z2) [fill=white] circle (\vr);
\draw (z3) [fill=white] circle (\vr);
\draw (z4) [fill=white] circle (\vr);
\draw[anchor = south] (u1) node {{\small $u_1$}};
\draw[anchor = south] (w2) node {{\small $v_1$}};
\draw[anchor = south] (x1) node {{\small $u_2$}};
\draw[anchor = south] (z2) node {{\small $v_2$}};
\end{tikzpicture}
\end{center}
\vskip -0.5 cm
\caption{An $\alpha$-set in the claw-free, diamond-free cubic graph $H_{2,\cub}$}
\label{fig:H3cubic}
\end{figure}

\begin{figure}[htb]
\begin{center}
\begin{tikzpicture}[scale=.8,style=thick,x=0.8cm,y=0.8cm]
\def\vr{2.5pt} 
\path (2.5,5.25) coordinate (v1);
\path (2.5,6.25) coordinate (v2);
\path (4,7.25) coordinate (v3);
\path (1,7.25) coordinate (v4);
\path (0,0.75) coordinate (u1);
\path (1,0) coordinate (u2);
\path (1,1.5) coordinate (u3);
\path (2,0.75) coordinate (u4);
\path (0,2.5) coordinate (u5);
\path (2,2.5) coordinate (u6);
\path (1,3.25) coordinate (u7);
\path (1,4.25) coordinate (u8);
\path (3,0.75) coordinate (w1);
\path (4,0) coordinate (w2);
\path (4,1.5) coordinate (w3);
\path (5,0.75) coordinate (w4);
\path (3,2.5) coordinate (w5);
\path (5,2.5) coordinate (w6);
\path (4,3.25) coordinate (w7);
\path (4,4.25) coordinate (w8);
\draw (u1)--(u2)--(u3)--(u4)--(u2)--(u3)--(u1);
\draw (u7)--(u5)--(u6)--(u7)--(u8);
\draw (u1)--(u5);
\draw (u4)--(u6);
\draw (w1)--(w2)--(w3)--(w4)--(w2)--(w3)--(w1);
\draw (w7)--(w5)--(w6)--(w7)--(w8);
\draw (w1)--(w5);
\draw (w4)--(w6);
\draw (v2)--(v3)--(v4)--(v2)--(v1);
\draw (v1)--(u8);
\draw (u8)--(w8);
\draw (v1)--(w8);
%
\draw (v1) [fill=white] circle (\vr);
\draw (v2) [fill=black] circle (\vr);
\draw (v3) [fill=white] circle (\vr);
\draw (v4) [fill=white] circle (\vr);
\draw (u1) [fill=white] circle (\vr);
\draw (u2) [fill=black] circle (\vr);
\draw (u3) [fill=white] circle (\vr);
\draw (u4) [fill=white] circle (\vr);
\draw (u5) [fill=white] circle (\vr);
\draw (u6) [fill=white] circle (\vr);
\draw (u7) [fill=black] circle (\vr);
\draw (u8) [fill=white] circle (\vr);
\draw (w1) [fill=white] circle (\vr);
\draw (w2) [fill=black] circle (\vr);
\draw (w3) [fill=white] circle (\vr);
\draw (w4) [fill=white] circle (\vr);
\draw (w5) [fill=white] circle (\vr);
\draw (w6) [fill=white] circle (\vr);
\draw (w7) [fill=black] circle (\vr);
\draw (w8) [fill=white] circle (\vr);
\draw[anchor = south] (v3) node {{\small $v$}};
\draw[anchor = south] (v4) node {{\small $u$}};
\draw (2.5,-1) node {{\small (a) A $\gamma$-set in $G_{20}$}};
\path (10.5,5.25) coordinate (v1);
\path (10.5,6.25) coordinate (v2);
\path (12,7.25) coordinate (v3);
\path (9,7.25) coordinate (v4);
\path (8,0.75) coordinate (u1);
\path (9,0) coordinate (u2);
\path (9,1.5) coordinate (u3);
\path (10,0.75) coordinate (u4);
\path (8,2.5) coordinate (u5);
\path (10,2.5) coordinate (u6);
\path (9,3.25) coordinate (u7);
\path (9,4.25) coordinate (u8);
\path (11,0.75) coordinate (w1);
\path (12,0) coordinate (w2);
\path (12,1.5) coordinate (w3);
\path (13,0.75) coordinate (w4);
\path (11,2.5) coordinate (w5);
\path (13,2.5) coordinate (w6);
\path (12,3.25) coordinate (w7);
\path (12,4.25) coordinate (w8);
\draw (u1)--(u2)--(u3)--(u4)--(u2)--(u3)--(u1);
\draw (u7)--(u5)--(u6)--(u7)--(u8);
\draw (u1)--(u5);
\draw (u4)--(u6);
\draw (w1)--(w2)--(w3)--(w4)--(w2)--(w3)--(w1);
\draw (w7)--(w5)--(w6)--(w7)--(w8);
\draw (w1)--(w5);
\draw (w4)--(w6);
\draw (v2)--(v3)--(v4)--(v2)--(v1);
\draw (v1)--(u8);
\draw (u8)--(w8);
\draw (v1)--(w8);
\draw (v1) [fill=black] circle (\vr);
\draw (v2) [fill=white] circle (\vr);
\draw (v3) [fill=white] circle (\vr);
\draw (v4) [fill=black] circle (\vr);
\draw (u1) [fill=black] circle (\vr);
\draw (u2) [fill=white] circle (\vr);
\draw (u3) [fill=white] circle (\vr);
\draw (u4) [fill=black] circle (\vr);
\draw (u5) [fill=white] circle (\vr);
\draw (u6) [fill=white] circle (\vr);
\draw (u7) [fill=black] circle (\vr);
\draw (u8) [fill=white] circle (\vr);
\draw (w1) [fill=black] circle (\vr);
\draw (w2) [fill=white] circle (\vr);
\draw (w3) [fill=white] circle (\vr);
\draw (w4) [fill=black] circle (\vr);
\draw (w5) [fill=white] circle (\vr);
\draw (w6) [fill=white] circle (\vr);
\draw (w7) [fill=black] circle (\vr);
\draw (w8) [fill=white] circle (\vr);
\draw[anchor = south] (v3) node {{\small $v$}};
\draw[anchor = south] (v4) node {{\small $u$}};
\draw (10.5,-1) node {{\small (b) An $\alpha$-set in $G_{20}$}};
\end{tikzpicture}
\end{center}
\vskip -0.5 cm
\caption{The claw-free, subcubic graph $G_{20}$}
\label{fig:G20}
\end{figure}
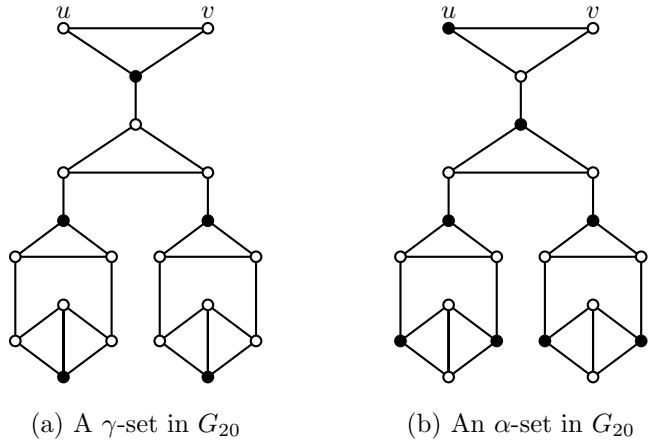

We now consider $k \ge 1$ copies of $G_{12}$, where the $i$th copy is denoted by $G_{12.i}$ and where $u_i$ and $v_i$ denote the two vertices of degree~$2$ in $G_{12.i}$ corresponding to the vertices $u$ and $v$ in $G_{12}$ indicated in Figure~\ref{fig:G20} for $i \in [k]$. Let $H_{1,\cub}$ be the graph obtained from $G_{12.1}$ by adding the edge~$u_1v_1$. For $k \ge 2$, let $H_{k,\cub}$ be the graph obtained from the disjoint union of these $k$ copies of $G_{12}$ by adding the edges $v_iu_{i+1}$ for $i \in [k-1]$ and adding the edge $v_kv_1$. Let $\cH_\cub = \{ H_{k,\cub} \colon k \ge 1\}$. For example, the graph $H_{2,\cub} \in \cH_\cub$ is illustrated in Figure~\ref{fig:H3cubic}. If $G \in \cH_\cub$, then $G = H_{k,\cub}$ for some $k \ge 1$ and $\gamma(G) = 3k$ and $\alpha(G) = 4k$, implying that $\alpha(G) = \frac{4}{3} \gamma(G)$. We therefore readily infer that equality holds in Theorem~\ref{thm:alpha-versus-dom-cubic} for every graph $G \in \cH_\cub$.

\begin{proposition}
\label{prop:3regularH}
If $G \in \cH_\cub$, then $G$ is a connected, claw-free, diamond-free, cubic graph satisfying $\alpha(G) = \frac{4}{3} \gamma(G)$.
\end{proposition}

Let $G_{20}$ be the claw-free, subcubic graph shown in Figure~\ref{fig:G20}. The graph $G_{20}$ satisfies $\gamma(G_{20}) = 5$ and $\alpha(G_{20}) = 8$, where the shaded vertices in Figure~\ref{fig:G20}(a) and~\ref{fig:G20}(b) are examples of a $\gamma$-set and an $\alpha$-set, respectively, of $G_{20}$.

We now consider $k \ge 2$ copies of $G_{20}$, where the $i$th copy is denoted by $G_{20.i}$ and where $u_i$ and $v_i$ denote the two vertices of degree~$2$ in $G_{20.i}$ corresponding to the vertices $u$ and $v$ in $G_{20}$ indicated in Figure~\ref{fig:G20}. Let $G_{k,\cub}$ be the graph obtained from the disjoint union of these $k$ copies of $G_{20}$ by adding the edges $v_iu_{i+1}$ for $i \in [k-1]$ and adding the edge $v_kv_1$. Let $\cG_\cub = \{ G_{k,\cub} \colon k \ge 2\}$. For example, the graph $G_{3,\cub} \in \cG_\cub$ is illustrated in Figure~\ref{fig:3cubic}. If $G \in \cG_\cub$, then $G = G_{k,\cub}$ for some $k \ge 2$ and $\gamma(G) = 5k$ and $\alpha(G) = 8k$, implying that $\alpha(G) = \frac{8}{5} \gamma(G)$. We therefore readily infer that equality holds in Theorem~\ref{thm:alpha-versus-dom} in the case when $d = 3$ for every graph $G \in \cG_\cub$.

\begin{proposition}
\label{prop:3regular}
If $G$ is a connected, claw-free, cubic graph of order~$n$, then $\alpha(G) \le \frac{8}{5} \gamma(G)$. Furthermore if $G \in \cG_\cub$, then $\alpha(G) = \frac{8}{5} \gamma(G)$.
\end{proposition}

\begin{figure}[htb]
\begin{center}
\begin{tikzpicture}[scale=.8,style=thick,x=0.8cm,y=0.8cm]
\def\vr{2.5pt} 
\path (2.5,5.25) coordinate (v1);
\path (2.5,6.25) coordinate (v2);
\path (4,7.25) coordinate (v3);
\path (1,7.25) coordinate (v4);
\path (1,8.25) coordinate (v5);
\path (18,7.25) coordinate (x4);
\path (18,8.25) coordinate (x5);
\path (0,0.75) coordinate (u1);
\path (1,0) coordinate (u2);
\path (1,1.5) coordinate (u3);
\path (2,0.75) coordinate (u4);
\path (0,2.5) coordinate (u5);
\path (2,2.5) coordinate (u6);
\path (1,3.25) coordinate (u7);
\path (1,4.25) coordinate (u8);
\path (3,0.75) coordinate (w1);
\path (4,0) coordinate (w2);
\path (4,1.5) coordinate (w3);
\path (5,0.75) coordinate (w4);
\path (3,2.5) coordinate (w5);
\path (5,2.5) coordinate (w6);
\path (4,3.25) coordinate (w7);
\path (4,4.25) coordinate (w8);
\draw (x4) to[out=0,in=0, distance=0.65cm] (x5);
\draw (v4) to[out=180,in=180, distance=0.65cm] (v5);
\draw (v5)--(x5);
\draw (u1)--(u2)--(u3)--(u4)--(u2)--(u3)--(u1);
\draw (u7)--(u5)--(u6)--(u7)--(u8);
\draw (u1)--(u5);
\draw (u4)--(u6);
\draw (w1)--(w2)--(w3)--(w4)--(w2)--(w3)--(w1);
\draw (w7)--(w5)--(w6)--(w7)--(w8);
\draw (w1)--(w5);
\draw (w4)--(w6);
\draw (v2)--(v3)--(v4)--(v2)--(v1);
\draw (v1)--(u8);
\draw (u8)--(w8);
\draw (v1)--(w8);
\draw (4,7.25)--(8,7.25);
\draw (11,7.25)--(15,7.25);
\draw (v1) [fill=black] circle (\vr);
\draw (v2) [fill=white] circle (\vr);
\draw (v3) [fill=white] circle (\vr);
\draw (v4) [fill=black] circle (\vr);
\draw (u1) [fill=black] circle (\vr);
\draw (u2) [fill=white] circle (\vr);
\draw (u3) [fill=white] circle (\vr);
\draw (u4) [fill=black] circle (\vr);
\draw (u5) [fill=white] circle (\vr);
\draw (u6) [fill=white] circle (\vr);
\draw (u7) [fill=black] circle (\vr);
\draw (u8) [fill=white] circle (\vr);
\draw (w1) [fill=black] circle (\vr);
\draw (w2) [fill=white] circle (\vr);
\draw (w3) [fill=white] circle (\vr);
\draw (w4) [fill=black] circle (\vr);
\draw (w5) [fill=white] circle (\vr);
\draw (w6) [fill=white] circle (\vr);
\draw (w7) [fill=black] circle (\vr);
\draw (w8) [fill=white] circle (\vr);
\draw[anchor = south] (v3) node {{\small $v_1$}};
\draw[anchor = south] (v4) node {{\small $u_1$}};
\path (9.5,5.25) coordinate (v1);
\path (9.5,6.25) coordinate (v2);
\path (11,7.25) coordinate (v3);
\path (8,7.25) coordinate (v4);
\path (7,0.75) coordinate (u1);
\path (8,0) coordinate (u2);
\path (8,1.5) coordinate (u3);
\path (9,0.75) coordinate (u4);
\path (7,2.5) coordinate (u5);
\path (9,2.5) coordinate (u6);
\path (8,3.25) coordinate (u7);
\path (8,4.25) coordinate (u8);
\path (10,0.75) coordinate (w1);
\path (11,0) coordinate (w2);
\path (11,1.5) coordinate (w3);
\path (12,0.75) coordinate (w4);
\path (10,2.5) coordinate (w5);
\path (12,2.5) coordinate (w6);
\path (11,3.25) coordinate (w7);
\path (11,4.25) coordinate (w8);
\draw (u1)--(u2)--(u3)--(u4)--(u2)--(u3)--(u1);
\draw (u7)--(u5)--(u6)--(u7)--(u8);
\draw (u1)--(u5);
\draw (u4)--(u6);
\draw (w1)--(w2)--(w3)--(w4)--(w2)--(w3)--(w1);
\draw (w7)--(w5)--(w6)--(w7)--(w8);
\draw (w1)--(w5);
\draw (w4)--(w6);
\draw (v2)--(v3)--(v4)--(v2)--(v1);
\draw (v1)--(u8);
\draw (u8)--(w8);
\draw (v1)--(w8);
\draw (v1) [fill=black] circle (\vr);
\draw (v2) [fill=white] circle (\vr);
\draw (v3) [fill=white] circle (\vr);
\draw (v4) [fill=black] circle (\vr);
\draw (u1) [fill=black] circle (\vr);
\draw (u2) [fill=white] circle (\vr);
\draw (u3) [fill=white] circle (\vr);
\draw (u4) [fill=black] circle (\vr);
\draw (u5) [fill=white] circle (\vr);
\draw (u6) [fill=white] circle (\vr);
\draw (u7) [fill=black] circle (\vr);
\draw (u8) [fill=white] circle (\vr);
\draw (w1) [fill=black] circle (\vr);
\draw (w2) [fill=white] circle (\vr);
\draw (w3) [fill=white] circle (\vr);
\draw (w4) [fill=black] circle (\vr);
\draw (w5) [fill=white] circle (\vr);
\draw (w6) [fill=white] circle (\vr);
\draw (w7) [fill=black] circle (\vr);
\draw (w8) [fill=white] circle (\vr);
\draw[anchor = south] (v3) node {{\small $v_2$}};
\draw[anchor = south] (v4) node {{\small $u_2$}};
\path (16.5,5.25) coordinate (v1);
\path (16.5,6.25) coordinate (v2);
\path (18,7.25) coordinate (v3);
\path (15,7.25) coordinate (v4);
\path (14,0.75) coordinate (u1);
\path (15,0) coordinate (u2);
\path (15,1.5) coordinate (u3);
\path (16,0.75) coordinate (u4);
\path (14,2.5) coordinate (u5);
\path (16,2.5) coordinate (u6);
\path (15,3.25) coordinate (u7);
\path (15,4.25) coordinate (u8);
\path (17,0.75) coordinate (w1);
\path (18,0) coordinate (w2);
\path (18,1.5) coordinate (w3);
\path (19,0.75) coordinate (w4);
\path (17,2.5) coordinate (w5);
\path (19,2.5) coordinate (w6);
\path (18,3.25) coordinate (w7);
\path (18,4.25) coordinate (w8);
\draw (u1)--(u2)--(u3)--(u4)--(u2)--(u3)--(u1);
\draw (u7)--(u5)--(u6)--(u7)--(u8);
\draw (u1)--(u5);
\draw (u4)--(u6);
\draw (w1)--(w2)--(w3)--(w4)--(w2)--(w3)--(w1);
\draw (w7)--(w5)--(w6)--(w7)--(w8);
\draw (w1)--(w5);
\draw (w4)--(w6);
\draw (v2)--(v3)--(v4)--(v2)--(v1);
\draw (v1)--(u8);
\draw (u8)--(w8);
\draw (v1)--(w8);
\draw (v1) [fill=black] circle (\vr);
\draw (v2) [fill=white] circle (\vr);
\draw (v3) [fill=white] circle (\vr);
\draw (v4) [fill=black] circle (\vr);
\draw (u1) [fill=black] circle (\vr);
\draw (u2) [fill=white] circle (\vr);
\draw (u3) [fill=white] circle (\vr);
\draw (u4) [fill=black] circle (\vr);
\draw (u5) [fill=white] circle (\vr);
\draw (u6) [fill=white] circle (\vr);
\draw (u7) [fill=black] circle (\vr);
\draw (u8) [fill=white] circle (\vr);
\draw (w1) [fill=black] circle (\vr);
\draw (w2) [fill=white] circle (\vr);
\draw (w3) [fill=white] circle (\vr);
\draw (w4) [fill=black] circle (\vr);
\draw (w5) [fill=white] circle (\vr);
\draw (w6) [fill=white] circle (\vr);
\draw (w7) [fill=black] circle (\vr);
\draw (w8) [fill=white] circle (\vr);
\draw[anchor = south] (v3) node {{\small $v_3$}};
\draw[anchor = south] (v4) node {{\small $u_3$}};
\end{tikzpicture}
\end{center}
\vskip -0.5 cm
\caption{An $\alpha$-set in the claw-free, cubic graph $G_{3,\cub} \in \cG_\cub$}
\label{fig:3cubic}
\end{figure}

\subsubsection{$4$-regular graphs}
\label{sec:4regular}

When $d = 4$, Theorem~\ref{thm:alpha-versus-dom} simplifies to $\alpha(G) \le \frac{5}{3} \gamma(G)$. We show that this bound is achievable.  Let $G = G_{15}$ be the claw-free, $4$-regular graph of order $n = 15$ shown in Figure~\ref{fig:G15}. By Observation~\ref{obs:bound-dom}, we infer that $\gamma(G) \ge \frac{1}{5}n = 3$. The shaded vertices in Figure~\ref{fig:G15}(a) form a dominating set of $G$ of cardinality~$3$, and hence $\gamma(G) = 3$. Moreover, by Theorem~\ref{thm:upper-bound-alpha}, we have $\alpha(G) \le \frac{1}{3}n = 5$. The shaded vertices in Figure~\ref{fig:G15}(b) form an independent set of $G$ of cardinality~$5$, and so $\alpha(G) = 5$. Hence, $\alpha(G) = \frac{5}{3}\gamma(G)$, that is, the graph $G$ achieves equality in Theorem~\ref{thm:alpha-versus-dom} in the case when $d = 4$.

\begin{figure}[htb]
\begin{center}
\begin{tikzpicture}[scale=.725,style=thick,x=0.725cm,y=0.725cm]
\def\vr{2.5pt} 
%
\path (.75, 2.5) coordinate (t5);
\path (-.75, 2.5) coordinate (t4);
\path (1.5, 1) coordinate (t3);
\path (-1.5, 1) coordinate (t2);
\path (0, -0.25) coordinate (t1);
\path (0, -2) coordinate (a1);
\path (-3.25, -3.25) coordinate (a2);
\path (3.25, -3.25) coordinate (a3);
\path (-2.75, 0.5) coordinate (b1);
\path (-5.25, 3.15) coordinate (b2);
\path (-1.5, 3.5) coordinate (c1);
\path (0, 6.5) coordinate (c2);
\path (1.5, 3.5) coordinate (d1);
\path (1.6, 3.5) coordinate (d1p);
\path (5.25, 3.15) coordinate (d2);
\path (2.75, 0.5) coordinate (e1);
\draw (t1)--(t4);
\draw (t1)--(t5);
\draw (t2)--(t3);
\draw (t2)--(t5);
\draw (t3)--(t4);
\draw (a1)--(t2);
\draw (a1)--(t3);
\draw (a1)--(a2);
\draw (a1)--(a3);
\draw (a2)--(a3);
\draw (b1)--(t1);
\draw (b1)--(t4);
\draw (b1)--(b2);
\draw (b1)--(a2);
\draw (b2)--(a2);
\draw (c1)--(b2);
\draw (c1)--(c2);
\draw (c1)--(t2);
\draw (c1)--(t5);
\draw (b2)--(c2);
\draw (d1)--(d2);
\draw (d1)--(c2);
\draw (d2)--(c2);
\draw (d1)--(t4);
\draw (d1)--(t3);
\draw (e1)--(d2);
\draw (e1)--(a3);
\draw (d2)--(a3);
\draw (e1)--(t5);
\draw (e1)--(t1);
\draw (t1) [fill=white] circle (\vr);
\draw (t2) [fill=white] circle (\vr);
\draw (t3) [fill=black] circle (\vr);
\draw (t4) [fill=white] circle (\vr);
\draw (t5) [fill=white] circle (\vr);
\draw (a1) [fill=white] circle (\vr);
\draw (a2) [fill=white] circle (\vr);
\draw (a3) [fill=white] circle (\vr);
\draw (b1) [fill=white] circle (\vr);
\draw (b2) [fill=black] circle (\vr);
\draw (c1) [fill=white] circle (\vr);
\draw (c2) [fill=white] circle (\vr);
\draw (d1) [fill=white] circle (\vr);
\draw (d2) [fill=white] circle (\vr);
\draw (e1) [fill=black] circle (\vr);
\draw (0.0,-4) node[anchor=north] {(a) A $\gamma$-set of $G_{15}$};
\draw[anchor = south] (c2) node {{\small $v$}};
\draw[anchor = south] (d1p) node {{\small $u$}};
\draw[anchor = west] (d2) node {{\small $w$}};
\begin{scope}[shift={(13,0)}]
\path (.75, 2.5) coordinate (t5);
\path (-.75, 2.5) coordinate (t4);
    \path (1.5, 1) coordinate (t3);
    \path (-1.5, 1) coordinate (t2);
    \path (0, -0.25) coordinate (t1);

    \path (0, -2) coordinate (a1);
    \path (-3.25, -3.25) coordinate (a2);
    \path (3.25, -3.25) coordinate (a3);

    \path (-2.75, 0.5) coordinate (b1);
    \path (-5.25, 3.15) coordinate (b2);

    \path (-1.5, 3.5) coordinate (c1);
    \path (0, 6.5) coordinate (c2);

    \path (1.5, 3.5) coordinate (d1);
    \path (1.6, 3.5) coordinate (d1p);
    \path (5.25, 3.15) coordinate (d2);

    \path (2.75, 0.5) coordinate (e1);

    \draw (t1)--(t4);
    \draw (t1)--(t5);
    \draw (t2)--(t3);
    \draw (t2)--(t5);
    \draw (t3)--(t4);

    \draw (a1)--(t2);
    \draw (a1)--(t3);
    \draw (a1)--(a2);
    \draw (a1)--(a3);
    \draw (a2)--(a3);

    \draw (b1)--(t1);
    \draw (b1)--(t4);
    \draw (b1)--(b2);
    \draw (b1)--(a2);
    \draw (b2)--(a2);

    \draw (c1)--(b2);
    \draw (c1)--(c2);
    \draw (c1)--(t2);
    \draw (c1)--(t5);
    \draw (b2)--(c2);

    \draw (d1)--(d2);
    \draw (d1)--(c2);
    \draw (d2)--(c2);
    \draw (d1)--(t4);
    \draw (d1)--(t3);

    \draw (e1)--(d2);
    \draw (e1)--(a3);
    \draw (d2)--(a3);
    \draw (e1)--(t5);
    \draw (e1)--(t1);

    \draw (t1) [fill=white] circle (\vr);
    \draw (t2) [fill=white] circle (\vr);
    \draw (t3) [fill=white] circle (\vr);
    \draw (t4) [fill=white] circle (\vr);
    \draw (t5) [fill=white] circle (\vr);

    \draw (a1) [fill=black] circle (\vr);
    \draw (a2) [fill=white] circle (\vr);
    \draw (a3) [fill=white] circle (\vr);

    \draw (b1) [fill=black] circle (\vr);
    \draw (b2) [fill=white] circle (\vr);

    \draw (c1) [fill=black] circle (\vr);
    \draw (c2) [fill=white] circle (\vr);

    \draw (d1) [fill=black] circle (\vr);
    \draw (d2) [fill=white] circle (\vr);

    \draw (e1) [fill=black] circle (\vr);
    \draw (0.0,-4) node[anchor=north] {(b)  A $\alpha$-set of $G_{15}$};
\draw[anchor = south] (c2) node {{\small $v$}};
\draw[anchor = south] (d1p) node {{\small $u$}};
\draw[anchor = west] (d2) node {{\small $w$}};
\end{scope}
\end{tikzpicture}
\end{center}
\vskip -0.75 cm
\caption{A $4$-regular, claw-free graph $G_{15}$}
\label{fig:G15}
\end{figure}

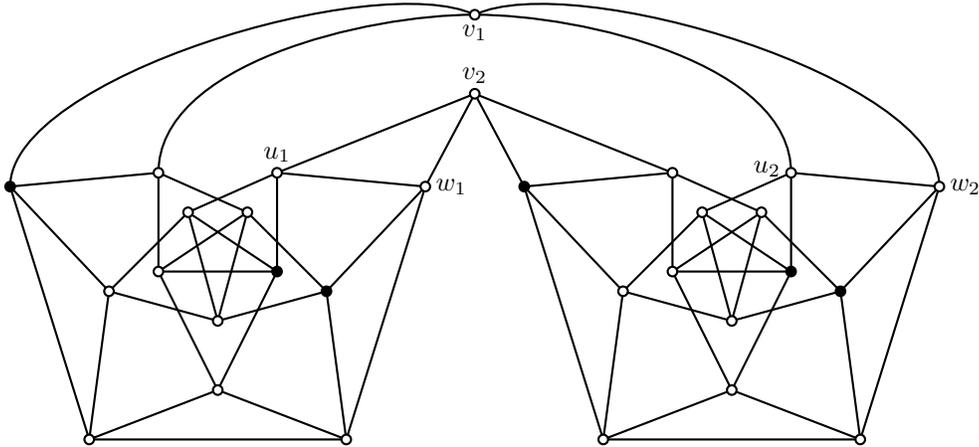
\begin{figure}[htb]
\begin{center}
\begin{tikzpicture}[scale=.725, style=thick, x=0.725cm, y=0.725cm]
\def\vr{2.5pt} 
\def\dx{13}     
\path ( 0.75,2.5) coordinate (t5_1);
\path (-0.75,2.5) coordinate (t4_1);
\path ( 1.5,1.0) coordinate (t3_1);
\path (-1.5,1.0) coordinate (t2_1);
\path (0.0,-0.25) coordinate (t1_1);
\path ( 0.0,-2.0)  coordinate (a1_1);
\path (-3.25,-3.25) coordinate (a2_1);
\path ( 3.25,-3.25) coordinate (a3_1);
\path (-2.75,0.5)  coordinate (b1_1);
\path (-5.25,3.15) coordinate (b2_1);
\path (-1.5,3.5)  coordinate (c1_1);
\path (6.5,7.5)  coordinate (c2_1);
\path (1.5,3.5)  coordinate (d1_1);
\path (5.25,3.15) coordinate (d2_1);
\path ( 2.75,  0.5)  coordinate (e1_1);
\path ( 0.75+\dx,  2.5)  coordinate (t5_2);
\path (-0.75+\dx,  2.5)  coordinate (t4_2);
\path ( 1.5+\dx,   1.0)  coordinate (t3_2);
\path (-1.5+\dx,   1.0)  coordinate (t2_2);
\path ( 0.0+\dx,  -0.25) coordinate (t1_2);
\path ( 0.0+\dx,  -2.0)  coordinate (a1_2);
\path (-3.25+\dx, -3.25) coordinate (a2_2);
\path ( 3.25+\dx, -3.25) coordinate (a3_2);
\path (-2.75+\dx,  0.5)  coordinate (b1_2);
\path (-5.25+\dx,  3.15) coordinate (b2_2);
\path (-1.5+\dx,   3.5)  coordinate (c1_2);
\path (-6.5+\dx,   5.5)  coordinate (c2_2);
\path ( 1.5+\dx,3.5)  coordinate (d1_2);
\path ( 1.5+\dx,3.6)  coordinate (d1_2p);
\path ( 5.25+\dx,3.15) coordinate (d2_2);
\path ( 2.75+\dx,  0.5)  coordinate (e1_2);
\draw (t1_1)--(t4_1);
\draw (t1_1)--(t5_1);
\draw (t2_1)--(t3_1);
\draw (t2_1)--(t5_1);
\draw (t3_1)--(t4_1);
\draw (a1_1)--(t2_1);
\draw (a1_1)--(t3_1);
\draw (a1_1)--(a2_1);
\draw (a1_1)--(a3_1);
\draw (a2_1)--(a3_1);
\draw (b1_1)--(t1_1);
\draw (b1_1)--(t4_1);
\draw (b1_1)--(b2_1);
\draw (b1_1)--(a2_1);
\draw (b2_1)--(a2_1);
\draw (c1_1)--(b2_1);
\draw (c1_1)--(t2_1);
\draw (c1_1)--(t5_1);
%
%
\draw (d1_1)--(d2_1);
\draw (d1_1)--(t4_1);
\draw (d1_1)--(t3_1);
\draw (e1_1)--(d2_1);
\draw (e1_1)--(a3_1);
\draw (d2_1)--(a3_1);
\draw (e1_1)--(t5_1);
\draw (e1_1)--(t1_1);
\draw[anchor = north] (c2_1) node {{\small $v_1$}};
\draw[anchor = south] (d1_1) node {{\small $u_1$}};
\draw[anchor = west] (d2_1) node {{\small $w_1$}};
\draw[anchor = south] (c2_2) node {{\small $v_2$}};
\draw[anchor = east] (d1_2p) node {{\small $u_2$}};
\draw[anchor = west] (d2_2) node {{\small $w_2$}};
\draw (c2_1) to[out=0,in=90, distance=2cm] (d1_2);
\draw (c2_1) to[out=25,in=90, distance=2cm] (d2_2);
\draw (c2_1) to[out=180,in=90, distance=2cm] (c1_1);
\draw (c2_1) to[out=155,in=90, distance=2cm] (b2_1);
\draw (t1_2)--(t4_2);
\draw (t1_2)--(t5_2);
\draw (t2_2)--(t3_2);
\draw (t2_2)--(t5_2);
\draw (t3_2)--(t4_2);

\draw (a1_2)--(t2_2);
\draw (a1_2)--(t3_2);
\draw (a1_2)--(a2_2);
\draw (a1_2)--(a3_2);
\draw (a2_2)--(a3_2);

\draw (b1_2)--(t1_2);
\draw (b1_2)--(t4_2);
\draw (b1_2)--(b2_2);
\draw (b1_2)--(a2_2);
\draw (b2_2)--(a2_2);

\draw (c1_2)--(b2_2);
\draw (c1_2)--(c2_2);
\draw (c1_2)--(t2_2);
\draw (c1_2)--(t5_2);
\draw (b2_2)--(c2_2);


\draw (d1_2)--(d2_2);
\draw (d1_2)--(t4_2);
\draw (d1_2)--(t3_2);

\draw (e1_2)--(d2_2);
\draw (e1_2)--(a3_2);
\draw (d2_2)--(a3_2);
\draw (e1_2)--(t5_2);
\draw (e1_2)--(t1_2);

\draw (c2_2)--(d1_1);
\draw (c2_2)--(d2_1);
%
\draw (t1_1) [fill=white] circle (\vr);
\draw (t2_1) [fill=white] circle (\vr);
\draw (t3_1) [fill=black] circle (\vr);
\draw (t4_1) [fill=white] circle (\vr);
\draw (t5_1) [fill=white] circle (\vr);
\draw (a1_1) [fill=white] circle (\vr);
\draw (a2_1) [fill=white] circle (\vr);
\draw (a3_1) [fill=white] circle (\vr);
\draw (b1_1) [fill=white] circle (\vr);
\draw (b2_1) [fill=black] circle (\vr);
\draw (c1_1) [fill=white] circle (\vr);
\draw (c2_1) [fill=white] circle (\vr);
\draw (d1_1) [fill=white] circle (\vr);
\draw (d2_1) [fill=white] circle (\vr);
\draw (e1_1) [fill=black] circle (\vr);
%
\draw (t1_2) [fill=white] circle (\vr);
\draw (t2_2) [fill=white] circle (\vr);
\draw (t3_2) [fill=black] circle (\vr);
\draw (t4_2) [fill=white] circle (\vr);
\draw (t5_2) [fill=white] circle (\vr);
\draw (a1_2) [fill=white] circle (\vr);
\draw (a2_2) [fill=white] circle (\vr);
\draw (a3_2) [fill=white] circle (\vr);
\draw (b1_2) [fill=white] circle (\vr);
\draw (b2_2) [fill=black] circle (\vr);
\draw (c1_2) [fill=white] circle (\vr);
\draw (c2_2) [fill=white] circle (\vr);
\draw (d1_2) [fill=white] circle (\vr);
\draw (d2_2) [fill=white] circle (\vr);
\draw (e1_2) [fill=black] circle (\vr);
\end{tikzpicture}
\end{center}
\vskip -0.25 cm
\caption{A $\gamma$-set in $G_{15}(2)$.}
\label{fig:G30bb}
\end{figure}

We now show how to construct arbitrarily large $4$-regular graphs that preserve the same ratio $\alpha(G) = \frac{5}{3}\gamma(G)$. Take $k \ge 2$ copies of $G_{15}$, denoted $G_1, G_2, \dots, G_k$, and label corresponding vertices in each copy as follows: denote by $v_j$ the ``top'' vertex in the $j$-th copy and by $u_j$ and $w_j$ the two vertices on the right that form a triangle with $v_j$ in $G_j$. Thus, the vertices $u_j$, $v_j$ and $w_j$ in $G_j$ correspond to the vertices named $u$, $v$ and $w$ in $G$ indicated in Figure~\ref{fig:G15}. For each copy $G_j$ where $j \in [k]$, we delete the edges $v_ju_j$ and $v_jw_j$ and we add the edges $v_{j+1}u_j$ and $v_{j+1}w_j$ cyclically, where we regard $v_{k+1}$ as $v_1$. This yields a new graph, which we denote $G_{15}(k)$. Observe that $G_{15}(k)$ has $15k$ vertices and is $4$-regular, connected, and claw-free (with each vertex neighborhood forming a ``constant link'' isomorphic to $2K_2$). Moreover, by construction each copy of $G_{15}$ still contributes an independent set of size~$5$, giving an overall independence number $\alpha(G_{15}(k)) = 5k$, and a dominating set of size~$3$ in each copy, giving $\gamma(G_{15}(k)) = 3k$. For example, a $\gamma$-set in $G_{15}(2)$ is indicated by the shaded vertices in Figure~\ref{fig:G30bb}. Thus we preserve the ratio
\[
\frac{\gamma(G_{15}(k))}{\alpha(G_{15}(k))} \;=\; \frac{3k}{5k} \;=\; \frac{3}{5}.
\]

\begin{proposition}
\label{prop:4regular-large}
There exist arbitrarily large connected, $4$-regular, claw-free graphs $G$ that satisfy $\alpha(G) = \frac{5}{3}\gamma(G)$.
\end{proposition}

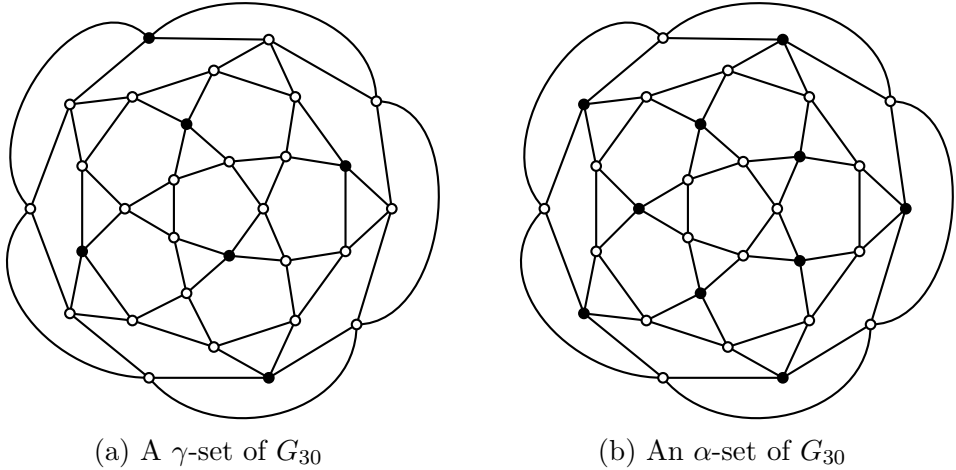
\begin{figure}[htb]
\begin{center}
\begin{tikzpicture}[scale=.725,style=thick,x=0.725cm,y=0.725cm]
\def\vr{2.5pt} 
%
\path (2,0) coordinate (t5);
\path (-1,4.28) coordinate (t4);
\path (2,8.6) coordinate (t3);
\path (7.75,7) coordinate (t2);
\path (7.25,1.35) coordinate (t1);
\path (5.03,0) coordinate (u1);
\path (8.14,4.28) coordinate (u2);
\path (5.03,8.56) coordinate (u3);
\path (0,6.92) coordinate (u4);
\path (0,1.63) coordinate (u5);
\path (3.64,0.78) coordinate (v1);
\path (5.70,1.45) coordinate (v2);
\path (6.97,3.20) coordinate (v3);
\path (6.97,5.36) coordinate (v4);
\path (5.70,7.11) coordinate (v5);
\path (3.64,7.78) coordinate (v6);
\path (1.58,7.11) coordinate (v7);
\path (0.31,5.36) coordinate (v8);
\path (0.31,3.20) coordinate (v9);
\path (1.58,1.45) coordinate (v10);
\path (2.95,2.14) coordinate (w1);
\path (5.46,2.96) coordinate (w2);
\path (5.46,5.60) coordinate (w3);
\path (2.95,6.42) coordinate (w4);
\path (1.39,4.28) coordinate (w5);
\path (4.03,3.09) coordinate (x1);
\path (4.89,4.28) coordinate (x2);
\path (4.03,5.47) coordinate (x3);
\path (2.63,5.01) coordinate (x4);
\path (2.63,3.55) coordinate (x5);
\draw (v1)--(v2)--(v3)--(v4)--(v5)--(v6)--(v7)--(v8)--(v9)--(v10)--(v1);
\draw (x1)--(x2)--(x3)--(x4)--(x5)--(x1);
\draw (v1)--(u1)--(v2);
\draw (v3)--(u2)--(v4);
\draw (v5)--(u3)--(v6);
\draw (v7)--(u4)--(v8);
\draw (v9)--(u5)--(v10);
\draw (v10)--(w1)--(v1);
\draw (v2)--(w2)--(v3);
\draw (v4)--(w3)--(v5);
\draw (v6)--(w4)--(v7);
\draw (v8)--(w5)--(v9);
\draw (w1)--(x1)--(w2)--(x2)--(w3)--(x3)--(w4)--(x4)--(w5)--(x5)--(w1);
\draw (u4)--(t4)--(u5);
\draw (u3)--(t3)--(u4);
\draw (u2)--(t2)--(u3);
\draw (u1)--(t1)--(u2);
\draw (u1)--(t5)--(u5);
\draw (t5) to[out=310,in=270, distance=1.75cm] (t1);
\draw (t1) to[out=0,in=0, distance=1.75cm] (t2);
\draw (t2) to[out=90,in=45, distance=1.75cm] (t3);
\draw (t3) to[out=135,in=135, distance=1.75cm] (t4);
\draw (t4) to[out=225,in=180, distance=1.75cm] (t5);
\draw (t1) [fill=white] circle (\vr);
\draw (t2) [fill=white] circle (\vr);
\draw (t3) [fill=black] circle (\vr);
\draw (t4) [fill=white] circle (\vr);
\draw (t5) [fill=white] circle (\vr);
\draw (u1) [fill=black] circle (\vr);
\draw (u2) [fill=white] circle (\vr);
\draw (u3) [fill=white] circle (\vr);
\draw (u4) [fill=white] circle (\vr);
\draw (u5) [fill=white] circle (\vr);
\draw (v1) [fill=white] circle (\vr);
\draw (v2) [fill=white] circle (\vr);
\draw (v3) [fill=white] circle (\vr);
\draw (v4) [fill=black] circle (\vr);
\draw (v5) [fill=white] circle (\vr);
\draw (v6) [fill=white] circle (\vr);
\draw (v7) [fill=white] circle (\vr);
\draw (v8) [fill=white] circle (\vr);
\draw (v9) [fill=black] circle (\vr);
\draw (v10) [fill=white] circle (\vr);
\draw (w1) [fill=white] circle (\vr);
\draw (w2) [fill=white] circle (\vr);
\draw (w3) [fill=white] circle (\vr);
\draw (w4) [fill=black] circle (\vr);
\draw (w5) [fill=white] circle (\vr);
\draw (x1) [fill=black] circle (\vr);
\draw (x2) [fill=white] circle (\vr);
\draw (x3) [fill=white] circle (\vr);
\draw (x4) [fill=white] circle (\vr);
\draw (x5) [fill=white] circle (\vr);
\draw (3.5,-1.25) node[anchor=north] {(a) A $\gamma$-set of $G_{30}$};
\begin{scope}[shift={(13,0)}]
    \path (2,0) coordinate (t5);
    \path (-1,4.28) coordinate (t4);
    \path (2,8.6) coordinate (t3);
    \path (7.75,7) coordinate (t2);
    \path (7.25,1.35) coordinate (t1);
    \path (5.03,0) coordinate (u1);
    \path (8.14,4.28) coordinate (u2);
    \path (5.03,8.56) coordinate (u3);
    \path (0,6.92) coordinate (u4);
    \path (0,1.63) coordinate (u5);
    \path (3.64,0.78) coordinate (v1);
    \path (5.70,1.45) coordinate (v2);
    \path (6.97,3.20) coordinate (v3);
    \path (6.97,5.36) coordinate (v4);
    \path (5.70,7.11) coordinate (v5);
    \path (3.64,7.78) coordinate (v6);
    \path (1.58,7.11) coordinate (v7);
    \path (0.31,5.36) coordinate (v8);
    \path (0.31,3.20) coordinate (v9);
    \path (1.58,1.45) coordinate (v10);
    \path (2.95,2.14) coordinate (w1);
    \path (5.46,2.96) coordinate (w2);
    \path (5.46,5.60) coordinate (w3);
    \path (2.95,6.42) coordinate (w4);
    \path (1.39,4.28) coordinate (w5);
    \path (4.03,3.09) coordinate (x1);
    \path (4.89,4.28) coordinate (x2);
    \path (4.03,5.47) coordinate (x3);
    \path (2.63,5.01) coordinate (x4);
    \path (2.63,3.55) coordinate (x5);
    \draw (v1)--(v2)--(v3)--(v4)--(v5)--(v6)--(v7)--(v8)--(v9)--(v10)--(v1);
    \draw (x1)--(x2)--(x3)--(x4)--(x5)--(x1);
    \draw (v1)--(u1)--(v2);
    \draw (v3)--(u2)--(v4);
    \draw (v5)--(u3)--(v6);
    \draw (v7)--(u4)--(v8);
    \draw (v9)--(u5)--(v10);
    \draw (v10)--(w1)--(v1);
    \draw (v2)--(w2)--(v3);
    \draw (v4)--(w3)--(v5);
    \draw (v6)--(w4)--(v7);
    \draw (v8)--(w5)--(v9);
    \draw (w1)--(x1)--(w2)--(x2)--(w3)--(x3)--(w4)--(x4)--(w5)--(x5)--(w1);
    \draw (u4)--(t4)--(u5);
    \draw (u3)--(t3)--(u4);
    \draw (u2)--(t2)--(u3);
    \draw (u1)--(t1)--(u2);
    \draw (u1)--(t5)--(u5);
    \draw (t5) to[out=310,in=270, distance=1.75cm] (t1);
    \draw (t1) to[out=0,in=0, distance=1.75cm] (t2);
    \draw (t2) to[out=90,in=45, distance=1.75cm] (t3);
    \draw (t3) to[out=135,in=135, distance=1.75cm] (t4);
    \draw (t4) to[out=225,in=180, distance=1.75cm] (t5);
    \draw (t1) [fill=white] circle (\vr);
    \draw (t2) [fill=white] circle (\vr);
    \draw (t3) [fill=white] circle (\vr);
    \draw (t4) [fill=white] circle (\vr);
    \draw (t5) [fill=white] circle (\vr);
    \draw (u1) [fill=black] circle (\vr);
    \draw (u2) [fill=black] circle (\vr);
    \draw (u3) [fill=black] circle (\vr);
    \draw (u4) [fill=black] circle (\vr);
    \draw (u5) [fill=black] circle (\vr);
    \draw (v1) [fill=white] circle (\vr);
    \draw (v2) [fill=white] circle (\vr);
    \draw (v3) [fill=white] circle (\vr);
    \draw (v4) [fill=white] circle (\vr);
    \draw (v5) [fill=white] circle (\vr);
    \draw (v6) [fill=white] circle (\vr);
    \draw (v7) [fill=white] circle (\vr);
    \draw (v8) [fill=white] circle (\vr);
    \draw (v9) [fill=white] circle (\vr);
    \draw (v10) [fill=white] circle (\vr);
    \draw (w1) [fill=black] circle (\vr);
    \draw (w2) [fill=black] circle (\vr);
    \draw (w3) [fill=black] circle (\vr);
    \draw (w4) [fill=black] circle (\vr);
    \draw (w5) [fill=black] circle (\vr);
    \draw (x1) [fill=white] circle (\vr);
    \draw (x2) [fill=white] circle (\vr);
    \draw (x3) [fill=white] circle (\vr);
    \draw (x4) [fill=white] circle (\vr);
    \draw (x5) [fill=white] circle (\vr);
\draw (3.5,-1.25) node[anchor=north] {(b)  An $\alpha$-set of $G_{30}$};
\end{scope}
\end{tikzpicture}
\end{center}
\vskip -0.75 cm
\caption{A $4$-regular, planar, claw-free graph $G_{30}$}
\label{fig:G30}
\end{figure}

The graph $G_{15}(k)$ is non-planar for any choice of $k$. Nevertheless, there exist $4$-regular, planar, claw-free graphs $G$ that satisfy $\alpha(G) = \frac{5}{3}\gamma(G)$, such as the graph $G_{30}$ illustrated in Figure~\ref{fig:G30}.



\section{On $k$-independence in $K_{1,r}$-free graphs}
\label{sec:k-independence}

In this section we provide a tight upper bound on the $k$-independence number of $K_{1,r}$-free graphs. Recall, a subset $S \subseteq V(G)$ is called a \emph{$k$-independent set} of a graph $G$ if the subgraph induced by $S$ has maximum degree at most $k$. The \emph{$k$-independence number} of $G$, denoted $\alpha_k(G)$, is the cardinality of a largest $k$-independent set in $G$. If $S \subseteq V(G)$ is a $k$-independent set with $|S| = \alpha_k(G)$, we refer to $S$ as an \emph{$\alpha_k$-set} of $G$. Note that $k$-independence generalizes the notion of independence: when $k = 0$, a $k$-independent set is simply an independent set, and thus $\alpha_0(G) = \alpha(G)$, the independence number of $G$. We next prove Theorem~\ref{thm:general_intro_kindep}, which provides a sharp upper bound for $\alpha_k(G)$ in $K_{1,r}$-free graphs in terms of the computable parameters $r$ and $k$, the minimum degree $\delta$, and the order $n$. Recall the statement of the theorem.

\smallskip
\noindent \textbf{Theorem~\ref{thm:general_intro_kindep}} \emph{For $k \ge 0$ and $r \ge 3$, if $G \in \cG_r$ has order $n$ and minimum degree $\delta \ge k+1$, then
\[
\alpha_k(G) \le \left( \frac{(r-1)(k+1)}{\delta - k + (r-1)(k+1)} \right) n,
\]
and this bound is sharp for all parameters involved, including for every $n \equiv 0 \pmod{(r-1)(k+1) + \delta - k}$.
}

\noindent
\proof For $k \ge 0$ and $r \ge 3$, let $G$ be a $K_{1,r}$-free graph of order $n$ and minimum degree $\delta \ge k+1$. Let $S \subseteq V(G)$ be an $\alpha_k$-set of $G$, and let $B = V(G) \setminus S$. Thus, $|S| = \alpha_k(G)$ and $|B| = n - \alpha_k(G)$. Let $e(S,B)$ be the number of edges with one end in $S$ and the other end in $B$. We note that each vertex $v \in S$ has at most~$k$ neighbors in~$S$ and at least $\deg_G(v) - k$ neighbors in $B$. Thus,
\[
e(S,B) \ge \sum_{v \in S}(\deg_G(v) - k)  \ge \sum_{v \in S}(\delta - k) = |S| (\delta - k) = \alpha_k(G)(\delta - k). \tag{4}
\]

We show next that each vertex in $B$ has at most~$(r-1)(k+1)$ neighbors in~$S$. Suppose, to the contrary, that there exists a vertex $v \in B$ that has at least $(r-1)(k+1) + 1$ neighbors in $S$, and let $H$ be the subgraph of $G$ induced by these neighbors of~$v$. Thus,
$|V(H)| \ge (r-1)(k+1) + 1$ and $\Delta(H) \le k$. Since $\chi(H) \le \Delta(H) + 1 \le k+1$ and $\chi(H)\alpha(H) \ge |V(H)|$, we infer that
\[
\alpha(H) \ge \left\lceil \frac{ |V(H)|}{ k+1 } \right\rceil = \left\lceil \frac{ (r-1)(k+1) + 1 }{ k+1 } \right\rceil = r.
\]

However if $S_v$ is an independent set in $H$ of cardinality~$r$, then the induced subgraph $G[S_v \cup \{v\}] = K_{1,r}$, contradicting the fact that $G$ is $K_{1,r}$-free. Hence each vertex in $B$ has at most~$(r-1)(k+1)$ neighbors in~$S$, implying that
\[
e(S,B) \le |B|(r-1)(k+1) = (r-1)(k+1)(n - \alpha_k(G)). \tag{5}
\]

Equating (4) and (5), we obtain $\alpha_k(G)(\delta - k) \le (r-1)(k+1)(n - \alpha_k(G))$. Solving for $\alpha_k(G)$, yields
\[
\alpha_k(G) \le \left( \frac{(r-1)(k+1)}{\delta - k + (r-1)(k+1)} \right) n,
\]
proving the desired bound.

To see that this bound is sharp, let $\delta = k + t$, where $t \ge 1$. Let $G$ be the graph obtained from the disjoint union of $r-1$ copies of $K_{k+1}$ by adding a complete graph $K_t$ and adding all edges between $K_t$ and the $r-1$ copies of $K_{k+1}$, that is, $G = (r-1)K_{k+1} + K_t$. Let $G$ have order~$n$, and so $n = (r-1)(k+1) + t$. In this construction of $G$, each copy of $K_{k+1}$, when combined with the vertices of $K_t$, forms a clique $K_{k+t+1}$. We note that every vertex in the added clique $K_t$ is a dominating vertex (of degree~$(r-1)(k+1) + t - 1$), while every vertex in one of the original copies of $K_{k+1}$ has degree~$k + t$. In particular, $\delta = k + t$ (where recall that $\delta = \delta(G)$). Every independent set in $G$ either contains a vertex in the added clique~$K_t$ (and no other vertex in $G$) or contains at most one vertex from each of the $r-1$ copies of $K_{k+1}$ in the construction of $G$, implying that $\alpha(G) \le r-1$. We therefore infer that the graph $G$ is $K_{1,r}$-free. The set $S$ consisting of the $(r-1)(k+1)$ vertices that belong to the $r-1$ copies of $K_{k+1}$ in the construction of $G$ is a $k$-independent set in $G$, and so
\[
\alpha_k(G) \ge (r-1)(k+1). \tag{6}
\]

As observed earlier, $\delta - k = t$ and $n = t + (r-1)(k+1)$. Hence, from the upper bound established in the theorem we have
\[
\begin{array}{lcl}
\alpha_k(G) & \le & \left( \frac{(r-1)(k+1)}{\delta - k + (r-1)(k+1)} \right) n \2  \\
& = & \left( \frac{(r-1)(k+1)}{t + (r-1)(k+1)} \right) \left( t + (r-1)(k+1) \right) \2  \\
& = &  (r-1)(k+1). \tag{7}
\end{array}
\]

Consequently, by (6) and (7), we infer that
\[
\alpha_k(G) = (r-1)(k+1) = \left( \frac{(r-1)(k+1)}{\delta - k + (r-1)(k+1)} \right) n.
\]

This construction demonstrates that the bound is sharp. Moreover, by taking multiple disjoint copies of the constructed graph, we see that the bound holds for infinitely many values of $n$, where $n \equiv 0 \pmod{(r-1)(k+1) + \delta - k}$.~\QED

\medskip
In view of the sharpness examples given above, it is of interest to present further constructions of arbitrarily large, connected graphs achieving the upper bound with the parameters given in the statement of Theorem~\ref{thm:general_intro_kindep}.

\begin{proposition}
\label{thm:yair-equality}
For $k \ge 0$ and $r \ge \max\{3,k+1\}$, there exist infinitely many connected $K_{1,r}$-free graphs of order~$n$ with minimum degree $\delta \ge k+2$ and $k$-independence number
\[
\alpha_k(G) = \left( \frac{(r-1)(k+1)}{\delta - k + (r-1)(k+1)} \right) n
\]
\end{proposition}
\proof
We begin by constructing the graph for the case $\delta = k+2$. Let $p$ be any positive integer. The graph consists of $p(r-1)$ disjoint copies of $K_{k+1}$, which are grouped into $p$ sets, each containing exactly $r-1$ copies of $K_{k+1}$. Let $S$ denote the set of all vertices that belong to these $p(r-1)$ copies of $K_{k+1}$. For each group indexed by $j \in [p]$, we introduce two vertices $x_{1,j}$ and $x_{2,j}$, and we define
\[
X = \bigcup_{j=1}^p \{x_{1,j},x_{2,j}\}.
\]

\begin{figure}[htb]
\begin{center}
\begin{tikzpicture}[scale=.85,style=thick,x=.85cm,y=.85cm]
\def\vr{2.5pt} 
\path (5.83,0.72) coordinate (v1);
\path (7.05,1.94) coordinate (v2);
\path (7.50,3.61) coordinate (v3);
\path (7.05,5.28) coordinate (v4);
\path (5.83,6.50) coordinate (v5);
\path (4.17,6.94) coordinate (v6);
\path (2.50,6.50) coordinate (v7);
\path (1.28,5.28) coordinate (v8);
\path (0.83,3.61) coordinate (v9);
\path (1.28,1.94) coordinate (v10);
\path (2.50,0.72) coordinate (v11);
\path (4.17,0.28) coordinate (v12);
\path (5.42,1.44) coordinate (u1);
\path (6.25,0) coordinate (w1);
\path (8.33,3.61) coordinate (u3);
\path (6.67,3.61) coordinate (w3);
\path (5.42,5.77) coordinate (u5);
\path (6.25,7.22) coordinate (w5);
\path (2.08,7.22) coordinate (u7);
\path (2.92,5.77) coordinate (w7);
\path (0,3.61) coordinate (u9);
\path (1.67,3.61) coordinate (w9);
\path (2.92,1.44) coordinate (u11);
\path (2.08,0) coordinate (w11);
\draw (u1)--(v1)--(w1);
\draw (u3)--(v3)--(w3);
\draw (u5)--(v5)--(w5);
\draw (u7)--(v7)--(w7);
\draw (u9)--(v9)--(w9);
\draw (u11)--(v11)--(w11);
\draw (v1)--(v2)--(v3)--(v4)--(v5)--(v6)--(v7)--(v8)--(v9)--(v10)--(v11)--(v12)--(v1);
\draw (w1) to[out=0,in=0, distance=0.35cm] (u1);
\draw (w3) to[out=90,in=90, distance=0.35cm] (u3);
\draw (w5) to[out=180,in=180, distance=0.35cm] (u5);
\draw (w7) to[out=180,in=180, distance=0.35cm] (u7);
\draw (w9) to[out=270,in=270, distance=0.35cm] (u9);
\draw (w11) to[out=0,in=0, distance=0.35cm] (u11);
\draw (u1)--(v2)--(w1);
\draw (u3)--(v2)--(w3);
\draw (u5)--(v4)--(w5);
\draw (u3)--(v4)--(w3);
\draw (u5)--(v6)--(w5);
\draw (u7)--(v6)--(w7);
\draw (u9)--(v8)--(w9);
\draw (u7)--(v8)--(w7);
\draw (u9)--(v10)--(w9);
\draw (u11)--(v10)--(w11);
\draw (u1)--(v12)--(w1);
\draw (u11)--(v12)--(w11);
\draw (v1) [fill=white] circle (\vr);
\draw (v2) [fill=black] circle (\vr);
\draw (v3) [fill=white] circle (\vr);
\draw (v4) [fill=black] circle (\vr);
\draw (v5) [fill=white] circle (\vr);
\draw (v6) [fill=black] circle (\vr);
\draw (v7) [fill=white] circle (\vr);
\draw (v8) [fill=black] circle (\vr);
\draw (v9) [fill=white] circle (\vr);
\draw (v10) [fill=black] circle (\vr);
\draw (v11) [fill=white] circle (\vr);
\draw (v12) [fill=black] circle (\vr);
\draw (u1) [fill=white] circle (\vr);
\draw (w1) [fill=white] circle (\vr);
\draw (u3) [fill=white] circle (\vr);
\draw (w3) [fill=white] circle (\vr);
\draw (u5) [fill=white] circle (\vr);
\draw (w5) [fill=white] circle (\vr);
\draw (u7) [fill=white] circle (\vr);
\draw (w7) [fill=white] circle (\vr);
\draw (u9) [fill=white] circle (\vr);
\draw (w9) [fill=white] circle (\vr);
\draw (u11) [fill=white] circle (\vr);
\draw (w11) [fill=white] circle (\vr);
\draw[anchor = west] (v2) node {{\small $x_{11}$}};
\draw[anchor = west] (v4) node {{\small $x_{21}$}};
\draw[anchor = south] (v6) node {{\small $x_{12}$}};
\draw[anchor = east] (v8) node {{\small $x_{22}$}};
\draw[anchor = east] (v10) node {{\small $x_{13}$}};
\draw[anchor = north] (v12) node {{\small $x_{23}$}};
\end{tikzpicture}
\end{center}
\vskip -0.35cm
\caption{An $\alpha$-set depicted in a graph $G$ appearing in the proof of Proposition~\ref{thm:yair-equality}. }
\label{fig:yair-elaborate-construction}
\end{figure}
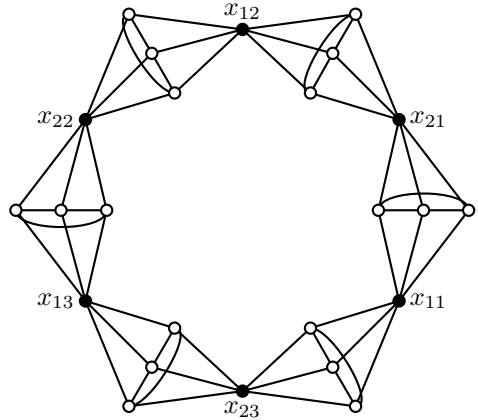

The vertex $x_{1,j}$ is joined to every vertex in the $j$-th group of $r-1$ copies of $K_{k+1}$, and the vertex $x_{2,j}$ is joined to every vertex within the $j$-th group of $r-1$ copies of $K_{k+1}$ except for the first copy of $K_{k+1}$ and to the first copy of $K_{k+1}$ in the $(j+1)$-th group (with indices considered cyclically). Finally, we add an edge between $x_{1,j}$ and $x_{2,j}$ for all $j \in [p]$. The resulting graph $G$ has order~$n = |X| + |S| = 2p + p(r-1)(k+1) = p\left(2 + (r-1)(k+1)\right)$. For example, when $r = 3$, $\delta = 4$, $k=2$, and $p=3$ the resulting graph $G$ is illustrated in Figure~\ref{fig:yair-elaborate-construction}.

By construction, the graph $G - X$ consists of $p(r-1)$ disjoint copies of $K_{k+1}$. Moreover every vertex in $G - X$ is adjacent to~$k$ neighbors in the copy of $K_{k+1}$ that contains it, and is adjacent to exactly two vertices in~$X$. Thus, every vertex in $G - X$ has degree~$k+2$ in $G$. Each vertex in $X$ has degree $(r-1)(k+1) + 1 > k + 2$. Thus, the graph $G$ has minimum degree $\delta = k+2$. By construction, each vertex in $X$ has at most~$r-1$ pairwise independent neighbors, implying that $G$ is $K_{1,r}$-free. Moreover, the set $S = V(G) \setminus X$ consisting of all vertices from the $p(r-1)$ copies of $K_{k+1}$ in the construction of $G$ forms a $k$-independent set in $G$, and so
\[
\alpha_k(G) \ge |S| = p(r-1)(k+1). \tag{8}
\]

Recall that $\delta - k = 2$. As observed earlier, $n = p\left(2 + (r-1)(k+1)\right)$. By Theorem~\ref{thm:general_intro_kindep}, we infer that
\[
\begin{array}{lcl}
\alpha_k(G) & \le & \left( \frac{(r-1)(k+1)}{\delta - k + (r-1)(k+1)} \right) n \2  \\
& = & \left( \frac{(r-1)(k+1)}{2 + (r-1)(k+1)} \right) \left( 2 + (r-1)(k+1)\right) p \2  \\
& = &  p(r-1)(k+1). \tag{9}
\end{array}
\]

Consequently, by (8) and (9), we infer that
\[
\alpha_k(G) = p(r-1)(k+1) = \left( \frac{(r-1)(k+1)}{\delta - k + (r-1)(k+1)} \right) n.
\]

Thus, equality holds in this case when $\delta = k + 2$. Now, we generalize the construction to $\delta = k + 2t + \mu$, where $\mu \in \{0, 1\}$ and $t \ge 1$. For each vertex $x_{1,j}$, we replace it with a clique of order $t + \mu$, and for each vertex $x_{2,j}$, we replace it with a clique of order $t$. These cliques are joined to the vertices of the $r-1$ copies of $K_{k+1}$ in the same way as $x_{1,j}$ and $x_{2,j}$ in the previous construction. Additionally, the edge between $x_{1,j}$ and $x_{2,j}$ is replaced by a clique of order $2t + \mu$.

In this generalized construction, each vertex in a $K_{k+1}$ component has degree at least $k + 2t + \mu$ because it is adjacent to $k$ vertices within its component and $2t + \mu$ additional vertices from the cliques replacing $x_{1,j}$ and $x_{2,j}$. Each vertex in the cliques replacing $x_{1,j}$ or $x_{2,j}$ has degree at least $(r-1)(k+1) + 2t - 1 + \mu$, which is strictly greater than $k + 2t + \mu$. Thus, the graph has minimum degree $\delta = k + 2t + \mu$ and remains $K_{1,r}$-free. The order of the generalized graph is
\[
n = p\left(2t + \mu + (r-1)(k+1)\right).
\]
Using the upper bound for $\alpha_k(G)$ given in Theorem~\ref{thm:general_intro_kindep} with $\delta - k = 2t + \mu$, we infer that
\[
\alpha_k(G) \le \left( \frac{(r-1)(k+1)}{\delta - k + (r-1)(k+1)} \right) n = p(r-1)(k+1).
\]
Equality again holds, completing the proof.~\QED

\section{Extremal bounds on induced subgraphs via Ramsey Theory}
\label{sec:ramsey}

The Ramsey-theoretic results in this section are all based on the same local observation. In a $K_{1,r}$-free graph, the neighborhood of every vertex has independence number at most~$r-1$. Thus, if a vertex is adjacent to sufficiently many vertices in an induced subgraph that avoids a prescribed family of subgraphs, Ramsey theory forces either one of the forbidden subgraphs inside that induced subgraph, or an independent set of cardinality~$r$ that forms an induced $K_{1,r}$ with the center vertex. We record this observation explicitly in the following lemma.

\begin{lemma}
\label{lem:local-ramsey-domination}
Let $r \ge 3$, let $G \in \cG_r$, and let $\mathcal{B}$ be a family of graphs. If $H$ is an induced subgraph of $G$ that contains no subgraph belonging to $\mathcal{B}$, and if $R = r(\mathcal{B},K_r)$, then the following hold.
\begin{enumerate}
\item If $x \in V(G) \setminus V(H)$, then $|N_G(x) \cap V(H)| \le R-1$.
\item If $x \in V(H)$, then $|N_H[x]| \le R-1$, and therefore $\deg_H(x) \le R-2$.
\end{enumerate}
\end{lemma}
\proof Suppose first that $x \in V(G) \setminus V(H)$ and that $X \subseteq N_G(x) \cap V(H)$ has cardinality~$R$. By the definition of $R$, the graph $G[X]$ contains either a subgraph in $\mathcal{B}$, impossible since $G[X]$ is a subgraph of $H$, or an independent set of cardinality~$r$, which together with $x$ induces a forbidden $K_{1,r}$. Thus $|N_G(x) \cap V(H)| \le R-1$.

Now let $x \in V(H)$ and suppose that $Y \subseteq N_H[x]$ has cardinality~$R$. Again, $G[Y]$ cannot contain a subgraph belonging to $\mathcal{B}$. Hence $G[Y]$ contains an independent set of cardinality~$r$. This independent set cannot contain $x$, since $x$ is adjacent to every other vertex in $N_H[x]$. It is therefore contained in $N_H(x)$ and, together with $x$, induces a forbidden $K_{1,r}$. Consequently $|N_H[x]| \le R-1$.~\QED

Before presenting the next result, we recall the notion of $\ell$-\emph{domination}. For an integer $\ell \ge 1$, a subset of vertices $D \subseteq V(G)$ is an \emph{$\ell$-dominating set} if every vertex outside $D$ has at least~$\ell$ neighbors in $D$. The minimum cardinality among all $\ell$-dominating sets in $G$ is the $\ell$-\emph{domination number} of $G$, denoted $\gamma_\ell(G)$. In particular, a $1$-dominating set is a dominating set of $G$, and $\gamma_1(G) = \gamma(G)$. An $\ell$-dominating set of $G$ with cardinality $\gamma_\ell(G)$ is called a $\gamma_\ell$-\emph{set of $G$}.

\begin{theorem}
\label{thm:2-domination}
For $\ell \in [3]$, if $G$ is a claw-free graph and $\mathcal{F}$ is the family of triangle-free graphs, then
\[
\alphaF(G) \le \frac{5}{\ell}\,\gamma_{\ell}(G),
\]
and this bound is sharp.
\end{theorem}
\proof Let $G$ be a claw-free graph and let $H$ be an induced subgraph of $G$ that is triangle-free and has maximum cardinality among all triangle-free induced subgraphs of $G$, and so $|V(H)| = \alphaF(G)$. For $\ell \in [3]$, let $D$ be a $\gamma_\ell$-set of $G$. Let
\[
S = D \cap V(H) \quad \mbox{and} \quad B = D \setminus S.
\]

Since $r(K_3,K_3)=6$, Lemma~\ref{lem:local-ramsey-domination}, applied with $\mathcal{B}=\{K_3\}$, implies that every vertex of $B$ is adjacent to at most five vertices of $H$. Further, since $H$ is triangle-free and $G$ is claw-free, every vertex of $H$ has degree at most~$2$ in $H$. Hence every vertex of $S$ is adjacent to at most two vertices of $V(H) \setminus S$. Let $e(D,V(H) \setminus S)$ denote the number of edges with one end in $D$ and the other in $V(H) \setminus S$. We have
\[
e(D,V(H) \setminus S) \le 5|B| + 2|S|.
\]
On the other hand, every vertex of $V(H) \setminus S$ has at least~$\ell$ neighbors in $D$, and therefore
\[
\ell|V(H) \setminus S| \le e(D,V(H) \setminus S).
\]
Since $\ell \le 3$, we obtain
\[
\begin{array}{lcl}
\ell|V(H)| & = & \ell|V(H) \setminus S| + \ell|S| \1 \\
& \le & 5|B| + 2|S| + \ell|S| \1 \\
& \le & 5(|B|+|S|) \1 \\
& = & 5|D|.
\end{array}
\]
Thus $\alphaF(G)=|V(H)| \le 5\gamma_\ell(G)/\ell$.

To see that this bound is sharp, let $G$ be obtained from a $5$-cycle $C_5$ by adding a copy of $K_\ell$, where $\ell \in [3]$, and adding all edges joining the added copy of $K_\ell$ and the $5$-cycle. The set of $\ell$ added vertices is a $\gamma_\ell$-set of $G$, and so $\gamma_\ell(G)=\ell$. Moreover, deleting these $\ell$ added vertices leaves the triangle-free cycle $C_5$, and no larger triangle-free induced subgraph exists by the bound just proved. Hence $\alphaF(G)=5=(5/\ell)\gamma_\ell(G)$.~\QED

\medskip
For $\ell \ge 4$, if $G$ is a claw-free graph and $\mathcal{F}$ is the family of triangle-free graphs, then it remains an open problem to determine a tight upper bound on $\alphaF(G)$ in terms of $\gamma_{\ell}(G)$. We remark that the proof ideas employed in the proof of Theorem~\ref{thm:2-domination} yield the bound $\alphaF(G) \le \left( \frac{\ell+2}{\ell} \right) \gamma_{\ell}(G)$, but it is unlikely that this bound is achievable.

Generalizing the result above, we now consider the case where $\mathcal{F}$ is the family of $K_q$-free graphs. We next prove Theorem~\ref{thm:general_intro_ramsey_kq}. Recall its statement.

\smallskip
\noindent \textbf{Theorem~\ref{thm:general_intro_ramsey_kq}} \emph{For integers $r, q \ge 3$, if $G \in \cG_{r}$ and if $\mathcal{F}$ is the family of $K_q$-free graphs, then
\[
\alphaF(G) \le \left(r(K_q,K_r) - 1\right)\gamma(G),
\]
and this bound is sharp.
}

\noindent
\proof Let $G$ be a $K_{1,r}$-free graph and let $\mathcal{F}$ be the family of $K_q$-free graphs. Let $H$ be a largest induced subgraph of $G$ that belongs to $\mathcal{F}$, and so $|V(H)|=\alphaF(G)$. Let $D$ be a $\gamma$-set of $G$, and put
\[
S = D \cap V(H) \quad \mbox{and} \quad B = D \setminus S.
\]
Set $R_q=r(K_q,K_r)$ and $R_{q-1}=r(K_{q-1},K_r)$.

By Lemma~\ref{lem:local-ramsey-domination}, applied with $\mathcal{B}=\{K_q\}$, every vertex of $B$ is adjacent to at most $R_q-1$ vertices of $H$. If $x \in S$, then $N_H(x)$ contains no $K_{q-1}$; otherwise such a copy of $K_{q-1}$ together with $x$ gives a copy of $K_q$ in $H$. Applying Lemma~\ref{lem:local-ramsey-domination} to this neighborhood gives $\deg_H(x) \le R_{q-1}-1$ for every $x \in S$.

Since $D$ dominates $G$, every vertex of $H$ lies either in the neighborhood of a vertex of $B$ or in the closed neighborhood in $H$ of a vertex of $S$. Hence
\[
\begin{array}{lcl}
|V(H)| & \le & (R_q-1)|B| + R_{q-1}|S| \1 \\
& \le & (R_q-1)(|B|+|S|) \1 \\
& = & (r(K_q,K_r)-1)\gamma(G),
\end{array}
\]
where we used the monotonicity $R_{q-1} \le R_q-1$ of Ramsey numbers. This proves the desired upper bound.

To see that this bound is sharp, consider a critical Ramsey graph $H$ for $r(K_q,K_r)$, and so $H$ has order $r(K_q,K_r)-1$, $H$ is $K_q$-free, and the complement $\barH$ is $K_r$-free. Equivalently, $H$ contains no independent set of cardinality~$r$. Let $G$ be the graph obtained from $H$ by adding a new vertex~$v$ adjacent to every vertex of $H$. The resulting graph $G$ is $K_{1,r}$-free and satisfies $\gamma(G)=1$. Since $H$ is an induced $K_q$-free subgraph of $G$, the upper bound just proved gives $\alphaF(G)=|V(H)|=r(K_q,K_r)-1=(r(K_q,K_r)-1)\gamma(G)$.~\QED

\medskip
Recall that a family of graphs $\mathcal{F}$ is \emph{edge}-\emph{hereditary} if for every graph $H \in \mathcal{F}$, every subgraph of $H$ also belongs to $\mathcal{F}$. Theorem~\ref{thm:general_intro_ramsey} is now an immediate consequence of Lemma~\ref{lem:local-ramsey-domination}. Recall its statement.

\smallskip
\noindent \textbf{Theorem~\ref{thm:general_intro_ramsey}} \emph{For $r \ge 3$, if $G \in \cG_{r}$ and $\mathcal{F}$ is an edge-hereditary family of graphs and $\mathcal{F}^*$ denotes the set of graphs not in $\mathcal{F}$, then
\[
\alphaF(G) \le \left(r(\mathcal{F}^*,K_r) - 1\right) \gamma(G),
\]
and this bound is sharp.
}

\noindent
\proof Let $G$ be a $K_{1,r}$-free graph and let $\mathcal{F}$ be an edge-hereditary family of graphs. Let $\mathcal{F}^*$ denote the set of graphs not in $\mathcal{F}$, and let $R=r(\mathcal{F}^*,K_r)$. Let $H$ be a largest induced subgraph of $G$ that belongs to $\mathcal{F}$, and so $|V(H)|=\alphaF(G)$. Since $\mathcal{F}$ is edge-hereditary, $H$ contains no subgraph that belongs to $\mathcal{F}^*$.

Let $D$ be a $\gamma$-set of $G$, and put
\[
S = D \cap V(H) \quad \mbox{and} \quad B = D \setminus S.
\]
By Lemma~\ref{lem:local-ramsey-domination}, every vertex of $B$ is adjacent to at most $R-1$ vertices of $H$, while every vertex of $S$ has closed neighborhood in $H$ of order at most $R-1$.

Since $D$ dominates $G$, every vertex of $H$ lies either in the neighborhood of a vertex of $B$ or in the closed neighborhood in $H$ of a vertex of $S$. Therefore
\[
\begin{array}{lcl}
|V(H)| & \le & (R-1)|B| + (R-1)|S| \1 \\
& = & (r(\mathcal{F}^*,K_r)-1)|D| \1 \\
& = & (r(\mathcal{F}^*,K_r)-1)\gamma(G).
\end{array}
\]
This proves the desired upper bound.

For sharpness, let $H$ be a critical Ramsey graph for $r(\mathcal{F}^*,K_r)$. Thus $H$ has order $r(\mathcal{F}^*,K_r)-1$, $H$ contains no subgraph from $\mathcal{F}^*$, and the complement $\barH$ is $K_r$-free; in particular, $H \in \mathcal{F}$ and $H$ contains no independent set of cardinality~$r$. Let $G$ be obtained from $H$ by adding a new vertex~$v$ adjacent to every vertex of $H$. Then $G$ is $K_{1,r}$-free and $\gamma(G)=1$. Since $H$ is an induced subgraph of $G$ that belongs to $\mathcal{F}$, the upper bound just proved gives $\alphaF(G)=|V(H)|=r(\mathcal{F}^*,K_r)-1=(r(\mathcal{F}^*,K_r)-1)\gamma(G)$.~\QED

\medskip
Our last result in this section derives a further relation between induced structures and domination via the same local Ramsey consideration. It is parallel in spirit to Theorem~\ref{thm:general_intro_kindep}, but with the domination parameter allowed to vary independently.

\begin{theorem}
\label{thm:kindep-ell-domination}
For $k \ge 0$, $\ell \ge 1$, and $r \ge 3$, if $G \in \cG_r$ and $\ell \le (r-1)(k+1)-k$, then
\[
\alpha_k(G) \le \frac{(r-1)(k+1)}{\ell}\,\gamma_\ell(G).
\]
\end{theorem}
\proof Let $H$ be a largest induced subgraph of $G$ with maximum degree at most~$k$, and so $|V(H)|=\alpha_k(G)$. Let $D$ be a $\gamma_\ell$-set of $G$, and put
\[
S = D \cap V(H) \quad \mbox{and} \quad B = D \setminus S.
\]
Set $C=(r-1)(k+1)$. Since $H$ contains no subgraph isomorphic to $K_{1,k+1}$, and since $r(K_{1,k+1},K_r)=C+1$ by Chv\'atal's theorem~\cite{Chvatal-1977}, Lemma~\ref{lem:local-ramsey-domination} implies that every vertex of $B$ is adjacent to at most $C$ vertices of $H$. Also, since $\Delta(H)\le k$, every vertex of $S$ is adjacent to at most $k$ vertices of $V(H)\setminus S$.

Every vertex of $V(H)\setminus S$ has at least~$\ell$ neighbors in $D$. Counting edges from $D$ to $V(H)\setminus S$, we obtain
\[
\ell|V(H)\setminus S| \le C|B|+k|S|.
\]
Thus,
\[
\begin{array}{lcl}
\ell|V(H)| & \le & C|B|+(\ell+k)|S| \1 \\
& \le & C(|B|+|S|) \1 \\
& = & C|D| \1 \\
& = & (r-1)(k+1)\gamma_\ell(G),
\end{array}
\]
where the second inequality follows from the hypothesis $\ell+k\le C$. The desired bound follows. ~\QED

\section{Conclusion}\label{sec:conclusion}

In this paper, we established sharp upper bounds on the independence number, $k$-independence number, and related invariants in $K_{1,r}$-free graphs in terms of the domination number. By introducing a general framework that bounds the maximum order of induced subgraphs with restricted chromatic number, we unified and extended several classical results. Notably, our results include tight bounds for bipartite, outerplanar, and planar subgraphs, as well as improvements in claw-free graphs and $r$-regular settings. Furthermore, we demonstrated that all bounds are sharp by providing explicit constructions achieving equality. However, several questions remain open.

\begin{problem}
\label{problem1}
{\rm For $d \ge 5$, determine a best possible upper bound on the independence number of a claw-free, $d$-regular graph $G$ in terms of its domination number. In particular, for $d \ge 5$, determine if the upper bound in Theorem~\ref{thm:alpha-versus-dom} is achievable.}
\end{problem}

\begin{problem}
\label{problem2}
{\rm For each value of $r \ge 3$, characterize the connected, $K_{1,r}$-free graphs $G$ satisfying $\alpha(G) = (r-1)\gamma(G)$.}
\end{problem}

\section*{Acknowledgements}

The authors thank the referees for their constructive suggestions, which helped simplify several proofs.

\section*{Statements and Declarations}

\section*{Funding}

The authors declare that no funding was received for this work.

\section*{Competing Interests}

The authors have no conflicts of interest to declare that are relevant to the content of this article.

\section*{Data availability}

No new data were created or analysed during this study. Data sharing is not applicable to this article.

\end{document}